\documentclass{article}\begin{document}

\centerline{\bf\large Riemann Spaces and Pfaff Differential Forms}
\vskip .4in

\centerline{\textbf{Nikos D. Bagis}}
\centerline{\textbf{Aristotele University of Thessaloniki-AUTH-Greece}}
\centerline{\textbf{nikosbagis@hotmail.gr}}
\[
\]
\centerline{\bf Abstract}
In this work we study differential geometry in $N$ dimensional Riemann curved spaces using Pfaff derivatives. Avoiding the classical partial derivative the Pfaff derivatives are constructed in a more sophisticated way and make evaluations become easier. In this way Christofell symbols $\Gamma_{ikj}$ of classical Riemann geometry as also the elements of the metric tensor $g_{ij}$ are replaced with one symbol (the $q_{ikj}$). Actually to describe the space we need no usage of the metric tensor $g_{ij}$ at all. We also don't use Einstein's notation and this quite simplifies things. For example we don't have to use upper and lower indexes, which in eyes of a beginner, is quite messy. Also we don't use the concept of tensor. All quantities of the surface or curve or space which form a tensor field are called invariants or curvatures of the space. Several new ideas are developed in this basis.

\[
\]

\section{Introduction and Development of the Theory}

Here we assume a $N$ dimensional space $\bf{S}$. The space $\bf{S}$ will be described by the vector
\begin{equation}
\overline{x}=\sum^{N}_{i=1}x_i(u_1,u_2,\ldots,u_N)\overline{\epsilon}_i ,
\end{equation}
where $\overline{\epsilon}_i$ is usual orthonormal base of $\textbf{E}=\textbf{R}^{N}$. We assume that in every point of the space $\bf{S}$ correspond $N$ orthonormal vectors $\{\overline{e}_1,\overline{e}_2,\ldots,\overline{e}_{N}\}$. 
These $N$ vectors $\{\overline{e}_1,\overline{e}_2,\ldots,\overline{e}_{N}\}$ span the space $\bf{S}$. We will use Pfaff derivatives to write our equations. We also study some properties of $\bf{S}$ which will need us for the construction of these equations. The Pfaff derivatives related with the structure of the space $\bf{S}$ which produces differential forms $\omega_k$, $k=1,2,\ldots,N$. These are defined as below.\\
It holds that 
\begin{equation}
\partial_j\overline{x}=\sum^{N}_{i=1}\partial_jx_i\overline{\epsilon}_i.
\end{equation}
The linear element of $\bf{S}$ is 
\begin{equation}
(ds)^2=(d\overline{x})^2=\sum^{N}_{i,j=1}\left\langle \partial_i\overline{x},\partial_j\overline{x}\right\rangle du_idu_j=\sum^{N}_{i=1}g_{ii}du_i^2+2\sum_{i<j}g_{ij}du_idu_j.
\end{equation} 
Hence
\begin{equation}
g_{ij}=\left\langle \frac{\partial \overline{x}}{\partial u_i},\frac{\partial \overline{x}}{\partial u_j} \right\rangle, 
\end{equation}
are the structure functions of the first linear form.\\ 
The Pfaff differential forms $\omega_k$ are defined with the help of $\{\overline{e}_k\}$, $k=1,2,\ldots,N$  as
\begin{equation}
d\overline{x}=\sum^{N}_{k=1}\omega_k\overline{e}_k.
\end{equation}

Then the Pfaff derivatives of the function $f$ are $\nabla_kf$, $k=1,2,\ldots,N$ and holds 
\begin{equation}
df=\sum^{N}_{k=1}(\partial_kf)du_k=\sum^{N}_{k=1}(\nabla_kf)\omega_k ,
\end{equation}
  
From (6) we get
\begin{equation}
d\overline{x}=\sum^{N}_{k=1}(\partial_k\overline{x})du_k.
\end{equation} 
Also 
\begin{equation}
d\overline{x}=\sum^{N}_{k=1}\omega_k\overline{e}_k\Rightarrow \nabla_m(\overline{x})=\overline{e}_m.
\end{equation}
Derivating the vectors $\overline{e}_l$ we can write them as linear composition of their self, since they form a basis of $\bf{E}$:
\begin{equation}
d\overline{e}_i=\sum^{N}_{k=1}\omega_{ik}\overline{e}_k.
\end{equation}
Then we define the connections $q_{ijm}$ and $b_{kl}$ as
\begin{equation}
\omega_{ij}=\sum^{N}_{m=1}q_{ijm}\omega_m\textrm{, }\omega_k=\sum^{N}_{l=1}b_{kl}du_{l}.
\end{equation}
Hence from (10),(5),(7) and $\left\langle \overline{e}_k ,d\overline{x}\right\rangle=\omega_k$, we get that
\begin{equation}
\left\langle \overline{e}_k, \partial_l\overline{x}\right\rangle=b_{kl}\textrm{ and }\partial_{l}\overline{x}=\sum^{N}_{k=1}b_{kl}\overline{e}_k
\end{equation}
and from (4)
\begin{equation}
g_{ij}=\sum^{N}_{k=1}b_{ki}b_{kj}. 
\end{equation}
Also
\begin{equation}
b_{sl}=\sum^{N}_{m=1}\frac{\partial x_m}{\partial u_l}\cos\left(\phi_{ms}\right)\textrm{, where }\cos\left(\phi_{ms}\right)=\left\langle \overline{\epsilon}_{m},\overline{e}_{s}\right\rangle
\end{equation}
and $\phi_{ms}$ is the angle formed by $\overline{\epsilon}_m$ and $\overline{e}_s$.
Also it holds
\begin{equation}
\omega_{ij}=\sum^{N}_{l=1}\left(\sum^{N}_{m=1}q_{ijm}b_{ml}\right)du_l
\end{equation} 
By this way we get the Christofell symbols  
\begin{equation}
\Gamma_{ijl}=\sum^{N}_{m=1}q_{ijm}b_{ml}.
\end{equation} 
Thus in view of (19) below 
\begin{equation}
\omega_{ij}=\sum^{N}_{l=1}\Gamma_{ijl}du_l\textrm{, }\Gamma_{ijl}+\Gamma_{jil}=0
\end{equation}
and easy (see and Proposition 1 below)
\begin{equation}
\partial_k\overline{e}_m=\sum^{N}_{j=1}\Gamma_{mjk}\overline{e}_j\textrm{ and }\nabla_k\overline{e}_m=\sum^{N}_{j=1}q_{mjk}\overline{e}_j.
\end{equation}

From the orthonormality of $\overline{e}_k$ we have  
\begin{equation}
\delta_{ij}=\left \langle \overline{e}_i,\overline{e}_j\right \rangle.
\end{equation} 
Derivating the above relation, we get
\begin{equation}
\omega_{ij}+\omega_{ji}=0
\end{equation}
and hence
\begin{equation}
q_{ijm}+q_{jim}=0.
\end{equation}
\\
\textbf{Theorem 1.}\\
The structure equations of $\bf{S}$ are (19) and 
\begin{equation}
d\omega_{j}=\sum^{N}_{m=1}\omega_{m}\wedge\omega_{mj}\textrm{ , }d\omega_{ij}=\sum^{N}_{m=1}\omega_{im}\wedge\omega_{mj}.
\end{equation}
\\
\textbf{Proof.}\\
We have
$$
d\left(d\overline{x}\right)=\overline{0}\Rightarrow d\left(\sum^{N}_{i=1}\overline{e}_i\omega_i\right)=\overline{0}\Rightarrow \sum^{N}_{i=1}\left(d\overline{e}_i\wedge\omega_i+\overline{e}_i d\omega_i\right)=\overline{0}.
$$
Hence
$$
\sum^{N}_{i=1}d\omega_i\overline{e}_i+\sum^{N}_{i,k=1}\omega_{ik}\wedge\omega_i\overline{e}_k=\overline{0}\Rightarrow d\omega_i=\sum^{N}_{l=1}\omega_{l}\wedge\omega_{li}.
$$
The same arguments hold and for the second relation of (21). $qed$\\
\\
\textbf{Definition 1.}\\
We write 
\begin{equation}
rot_{ij}\left(A_{k\ldots i\ldots j\ldots l}\right)=A_{k\ldots i\ldots j\ldots l}-A_{k\ldots j\ldots i\ldots l}
\end{equation}
\\
\textbf{Definition 2.}\\
We define the Kronecker$-\delta$ symbol as follows:\\
If $\{i_1,i_2,\ldots,i_M\}$, $\{j_1,j_2,\ldots,j_M\}$ are two set of indexes, then
$$
\delta_{{i_1i_2\ldots i_M},{j_1j_2\ldots j_M}}=1,
$$
if $\{i_1,i_2,\ldots,i_M\}$ is even permutation of  $\{j_1,j_2,\ldots,j_M\}$.
$$
\delta_{{i_1i_2\ldots i_M},{j_1j_2\ldots j_M}}=-1,
$$
if $\{i_1,i_2,\ldots,i_M\}$ is odd permutation of  $\{j_1,j_2,\ldots j_{M}\}$
and 
$$
\delta_{{i_1i_2\ldots i_M},{j_1j_2\ldots j_M}}=0
$$
in any other case.\\
\\

Let $a=\sum^{N}_{i=1}a_i\omega_i$ and $b=\sum^{N}_{j,k=1}b_{jk}\omega_j\wedge\omega_k$, then 
$$
(a\wedge b)_{123}=\frac{1}{1!}\frac{1}{2!}\sum^{N}_{i,j,k=1}\delta_{{ijk},{123}}a_jb_{jk}=
$$
$$
=\frac{1}{2}[a_1b_{23}\delta_{{123},{123}}+a_1b_{32}\delta_{{132},{123}}+a_2b_{13}\delta_{{213},{123}}+a_2b_{31}\delta_{{231},{123}}+
$$
$$
+a_3b_{12}\delta_{{312},{123}}+a_3b_{21}\delta_{{321},{123}}]=
$$
$$
=\frac{1}{2}\left(a_1b_{23}-a_1b_{32}-a_2b_{13}+a_2b_{31}+a_3b_{12}-a_3b_{21}\right).
$$ 
We remark here that we don't use Einstein's index notation.\\   
\\

By this way equations (21) give
$$
d\omega_j=\sum^{N}_{m=1}\omega_m\wedge\omega_{mj}=\sum^{N}_{m,s=1}q_{mjs}\omega_{m}\wedge\omega_{s}.
$$
Hence
\begin{equation}
d\omega_j=\sum_{m<s}Q_{mjs}\omega_m\wedge\omega_s,
\end{equation}
with 
\begin{equation}
Q_{mjs}:=rot_{ms}(q_{mjs})=q_{mjs}-q_{sjm} 
\end{equation}
Hence with the above notation holds $Q_{ikl}+Q_{lki}=0$. Also if we set
\begin{equation}
R_{ijkl}:=-\sum^{N}_{m=1}rot_{kl}\left(q_{imk}q_{jml}\right),
\end{equation}
then we have
\begin{equation}
d\omega_{ij}=\sum_{k<l}R_{ijkl}\omega_k\wedge\omega_l.
\end{equation}
\\
Using the structure equations of the space $\bf{S}$ we have the next\\
\\
\textbf{Proposition 1.}\\
\begin{equation}
\nabla_k\nabla_m(\overline{x})=\nabla_k(\overline{e}_m)=\sum^{N}_{j=1}q_{mjk}\overline{e}_j.
\end{equation} 
\textbf{Proof.}\\
$$
d\overline{e}_i=\sum^{N}_{k=1}\nabla_k(\overline{e}_i)\omega_k\Leftrightarrow  \sum^{N}_{k=1}\omega_{ik}\overline{e}_k=\sum^{N}_{k=1}\nabla_k(\overline{e}_i)\omega_k\Rightarrow
$$
$$ \sum^{N}_{m=1}\left(\sum^{N}_{k=1}q_{ikm}\overline{e}_k\right)\omega_m=\sum^{N}_{m=1}\nabla_{m}(\overline{e}_i)\omega_m.
$$ 
\\
\textbf{Theorem 2.}\\
For every function $f$ hold the following relations
\begin{equation}
\nabla_l\nabla_mf-\nabla_m\nabla_lf+\sum^{N}_{k=1}(\nabla_kf)Q_{lkm}=0\textrm{, }\forall\mbox{ } l,m\in\{1,2,\ldots,N\} 
\end{equation}
or equivalently
\begin{equation}
rot_{lm}\left(\nabla_l\nabla_mf+\sum^{N}_{k=1}(\nabla_kf)q_{lkm}\right)=0.
\end{equation}
\\
\textbf{Proof.}
$$
df=\sum^{N}_{k=1}(\nabla_kf)\omega_k\Rightarrow d(df)=\sum^{N}_{k=1}d(\nabla_kf)\wedge\omega_k+\sum^{N}_{k=1}(\nabla_kf)d\omega_k=0
$$
or
$$
\sum^{N}_{k=1}\left(\sum^{N}_{s=1}(\nabla_s\nabla_kf)\omega_s\right)\wedge\omega_k+\sum^{N}_{k=1}(\nabla_kf)\left(\sum_{l<m}rot_{lm}(q_{l km})\omega_{l}\wedge\omega_{m}\right)=0
$$
or
$$
\sum_{l<m}\left(\nabla_l\nabla_mf-\nabla_m\nabla_lf+\sum^{N}_{k=1}(\nabla_kf)Q_{lkm}\right)\omega_l\wedge\omega_m=0.
$$
\\
\textbf{Corollary 1.}\\
If $\lambda_{ij}=-\lambda_{ji}$, $i,j\in\{1,2,\ldots,N\}$ is any antisymetric field, then
\begin{equation}
\sum^{N}_{i,j=1}\lambda_{ij}\nabla_i\nabla_jf+\sum^{N}_{i,j,k=1}\lambda_{ij}q_{ikj}\nabla_kf=0.
\end{equation} 
\\
\textbf{Theorem 3.}
\begin{equation}
R_{iklm}=\nabla_l(q_{ikm})-\nabla_{m}(q_{ikl})+\sum^{N}_{s=1}q_{iks}Q_{lsm}=-\sum^{N}_{s=1}rot_{lm}(q_{isl}q_{ksm})
\end{equation}
and it holds also
\begin{equation}
R_{iklm}=-R_{kilm}\textrm{, }R_{iklm}=-R_{ikml}\textrm{, }
\end{equation}
\\
\textbf{Proof.}\\
We describe the proof.
$$
d\omega_{ij}=\sum^{N}_{k=1}d(q_{ijk})\wedge\omega_k+\sum^{N}_{k=1}q_{ijk}d\omega_k=
$$
$$
=\sum^{N}_{m<l}\left(\nabla_{m}(q_{ijl})-\nabla_{l}(q_{ijm})+\sum^{N}_{s=1}q_{ijs}Q_{msl}\right)\omega_m\wedge\omega_{l}
$$
and
$$
d\omega_{ij}=\sum^{N}_{k=1}\omega_{ik}\wedge\omega_{kj}
=\sum^{N}_{m<l}rot_{ml}\left(\sum^{N}_{s=1}q_{ism}q_{sjl}\right)\omega_{m}\wedge\omega_{l}
$$
Having in mind the above two relations we get the result.\\
\\
\textbf{Note 1.}\\
When exist a field $\Phi_{ij}$ such that $\nabla_{k}\Phi_{ij}=q_{ijk}$ : (a), then from Theorem 2 we have $d(\omega_{ij})=0$ and using Theorem 3: 
\begin{equation}
R_{iklm}=-\sum_{s=1}^{N}rot_{lm}\left(q_{isl}q_{ksm}\right)=0
\end{equation}
and in view of (26),(16):
\begin{equation}
\partial_l\Gamma_{ijk}-\partial_k\Gamma_{ijl}=0.
\end{equation}
\\
\textbf{Proposition 2.}
\begin{equation}
\nabla_l\nabla_m(\overline{e}_i)=\sum^{N}_{k=1}\left(\nabla_l\left(q_{ikm}\right)-\sum^{N}_{s=1}q_{ism}q_{ksl}\right)\overline{e}_k
\end{equation}
Hence
\begin{equation}
\left\langle \nabla_k^2(\overline{e}_m),\overline{e}_m \right\rangle=-\sum^{N}_{s=1}q_{msk}^2=\textrm{invariant}
\end{equation}
\begin{equation}
T_{ijlm}:=\left\langle \nabla_{l}\nabla_m(\overline{e}_i)-\nabla_m\nabla_{l}(\overline{e}_i),\overline{e}_{j}\right\rangle=\textrm{invariant}, 
\end{equation}
and
\begin{equation}
T_{ijlm}=rot_{lm}\left(\nabla_l\left(q_{ijm}\right)+\sum^{N}_{s=1}q_{isl}q_{jsm}\right)
\end{equation}
for all $i,j,l,m\in\{1,2,\ldots,N\}$.\\
\\
\textbf{Proof.}\\
See Lemma 1 below.\\
\\
\textbf{Theorem 4.}
\begin{equation}
T_{ijlm}=-\sum^{N}_{s=1}q_{ijs}Q_{lsm}
\end{equation}
and
\begin{equation}
\nabla_l\left(q_{ijm}\right)-\nabla_m\left(q_{ijl}\right)
\end{equation}
are invariants.\\
\\
\textbf{Proof.}\\
Use Theorem 3 along with Proposition 2.\\
\\
\textbf{Definition 3.}\\
We construct the differential operator 
$\Theta^{(1)}_{lm}$, such that for a vector field $\overline{Y}=Y_1\overline{e}_{1}+Y_2\overline{e}_2+\ldots+Y_N\overline{e}_N$ it is
\begin{equation}
\Theta^{(1)}_{lm}(\overline{Y}):=\nabla_l(Y_m)-\nabla_m(Y_l)+\sum^{N}_{k=1}Y_kQ_{lkm}.
\end{equation}
\\
\textbf{Theorem 5.}\\
The derivative of a vector  $\overline{Y}=Y_1\overline{e}_{1}+Y_2\overline{e}_2+\ldots+Y_N\overline{e}_N$ is
 \begin{equation}
d\overline{Y}=\sum^{N}_{j=1}\left(\sum^{N}_{l=1}Y_{j;l}\omega_l\right)\overline{e}_j 
\end{equation} 
where
\begin{equation}
Y_{j;l}=\nabla_lY_j-\sum^{N}_{k=1}q_{jkl}Y_k=\textrm{invariant}.
\end{equation}
\\
\textbf{Remark 1.}
\begin{equation}
\Theta_{lm}^{(1)}\left(\overline{Y}\right)=Y_{m;l}-Y_{l;m}
\end{equation}
\\
\textbf{Lemma 1.}\\
For every vector field $\overline{Y}$ we have
\begin{equation}
\left\langle \nabla_k \overline{Y},\overline{e}_l\right\rangle=Y_{l;k}=\textrm{invariant}
\end{equation}
\\
\textbf{Proof.}
$$
\left\langle \nabla_k \overline{Y},\overline{e}_l\right\rangle=\left\langle \nabla_k\left(\sum^{N}_{m=1}Y_m\overline{e}_m\right),\overline{e}_l\right\rangle=
$$
$$
=\sum^{N}_{m=1}\left\langle\nabla_k\left(Y_m\right)\overline{e}_m+Y_m\nabla_k\left(\overline{e}_m\right),\overline{e}_l\right\rangle=
$$
$$
=\nabla_kY_l+\sum^{N}_{m=1}Y_m\sum^{N}_{s=1}q_{msk}\left\langle\overline{e}_s,\overline{e}_l\right\rangle=\nabla_kY_l+\sum^{N}_{m=1}Y_mq_{mlk}=
$$
$$
=\nabla_kY_l-\sum^{N}_{m=1}Y_mq_{lmk}=Y_{l;k}=\textrm{invariant}.
$$
\\
\textbf{Proposition 3.}\\
The connections $q_{ijk}$ are invariants.\\
\\
\textbf{Proof.}
$$
\left\langle \nabla_k\overline{e}_m,\overline{e}_l\right\rangle=\left\langle\sum^{N}_{j=1}q_{mjk}\overline{e}_j,\overline{e}_l\right\rangle=q_{mlk}.
$$
\\
\textbf{Definition 4.}\\
If $\omega=\sum^{N}_{k=1}a_k\omega_k$ is a Pfaff form, then we define
\begin{equation}
\Theta^{(2)}_{lm}(\omega):=\nabla_{l}a_m-\nabla_{m}a_l+\sum^{N}_{k=1}a_kQ_{lkm}
\end{equation}
Hence 
\begin{equation}
d\omega=\sum^{N}_{l<m}\Theta^{(2)}_{lm}\left(\omega\right)\omega_{l}\wedge\omega_{m},
\end{equation}

The above derivative operators $X_{m;l}$, $\Theta^{(1)}_{lm}\left(\overline{Y}\right)$ and $\Theta^{(2)}_{lm}(\omega)$ are invariant under all acceptable change of variables. As application of such differentiation are the forms $\omega_{ij}=\sum^{N}_{k=1}q_{ijk}\omega_k$, which we can write
\begin{equation}
\Theta^{(2)}_{lm}\left(\omega_{ij}\right)=R_{ijlm}
\end{equation}
and lead us concluding that $R_{ijlm}$ are invariants of $\textbf{S}$. Actually\\
\\ 
\textbf{Theorem 6.}\\
$R_{ijlm}$ is the curvature tensor of $\textbf{S}$.\\   
\\
\textbf{Theorem 7.}\\
If $\omega=\sum^{N}_{k=1}a_{k}\omega_k$, then
\begin{equation}
d(f\omega)=\sum_{l<m}\left(\left|
\begin{array}{cc}
\nabla_lf\textrm{ }\nabla_{m}f\\
a_l\textrm{ }\textrm{ }a_m	
\end{array}
\right|+\Theta^{(2)}_{lm}(\omega)f\right)\omega_{l}\wedge\omega_{m}.
\end{equation}
In particular
\begin{equation}
d(f\omega_{ij})=\sum_{l<m}\left(\left|
\begin{array}{cc}
\nabla_lf\textrm{ }\nabla_{m}f\\
q_{ijl}\textrm{ }\textrm{ }q_{ijm}	
\end{array}
\right|+R_{ijlm}f\right)\omega_{l}\wedge\omega_{m}.
\end{equation}
Also
\begin{equation}
\Theta^{(2)}_{lm}(f\omega)=\left|
\begin{array}{cc}
	\nabla_lf\textrm{  }\nabla_mf\\
	a_l\textrm{  }\textrm{  }a_m
\end{array}
\right|+f\Theta^{(2)}_{lm}(\omega).
\end{equation}
If we set
\begin{equation}
R:=R_1\omega_1+R_2\omega_2+\ldots+R_N\omega_N,
\end{equation}
then
\begin{equation}
\Theta^{(2)}_{lm}\left(R\right)=\nabla_lR_m-\nabla_mR_l+\sum^{N}_{k=1}R_kQ_{lkm}
\end{equation}
and
\begin{equation}
\sum_{l<m}\Theta^{(2)}_{lm}(\omega)=\sum_{l<m}\left(\nabla_la_m-\nabla_ma_l\right)+\sum^{N}_{k=1}a_kA_k,
\end{equation}
where the $A_k$ are defined in Deffinition 6 bellow.\\
\\
\textbf{Theorem 8.}\\
We set 
\begin{equation}
\overline{R}:=R_1\overline{e}_1+R_2\overline{e}_2+\ldots+R_N\overline{e}_N,
\end{equation}
where the $R_k$ are as in Definition 6 below, then
\begin{equation}
\sum_{l<m}\left(\nabla_lR_m-\nabla_mR_l\right)=\textrm{invariant}.
\end{equation}
\\
\textbf{Proof.}\\
From Definition 6 below and $\sum^{N}_{k=1}R_kA_k=0$ (relation (76) below), we get
\begin{equation}
\sum_{l<m}\Theta^{(1)}_{lm}\left(\overline{R}\right)=\sum_{l<m}\left(\nabla_lR_m-\nabla_mR_l+\sum^{N}_{k=1}R_kQ_{lkm}\right)
\end{equation}
\\
\textbf{Note 2.}\\
\textbf{i)}
\begin{equation}
\Theta^{(2)}_{lm}(\omega_j)=Q_{ljm}=\textrm{invariant}
\end{equation}
Hence 
\begin{equation}
\sum_{l<m}\Theta^{(2)}_{lm}(\omega_j)=A_j=\textrm{invariant}
\end{equation}
\textbf{ii)}
If we assume that $\omega=dg=\sum^{N}_{k=1}(\nabla_{k}g)\omega_k$ and use (27) we get 
\begin{equation}
\Theta^{(2)}_{lm}(dg)=0
\end{equation}
and hence for all multivariable functions $f,g$ we have the next\\
\\
\textbf{Theorem 9.}
\begin{equation} 
\int_{\partial A}fdg=\sum_{l<m}\int\int_A\left(\nabla_lf\nabla_mg-\nabla_mf\nabla_lg\right)\omega_l\wedge\omega_m
\end{equation}
\\
\textbf{Proposition 4.}\\
If exists multivariable function $f=f(u_1,u_2,\ldots,u_N)\in\textbf{R}$ such that
\begin{equation}
\left|
\begin{array}{cc}
	\nabla_1f\textrm{ }\nabla_2f\\
	q_{ij1}\textrm{ }q_{ij2}
\end{array}
\right|=\left|
\begin{array}{cc}
	\nabla_2f\textrm{ }\nabla_3f\\
	q_{ij2}\textrm{ }q_{ij3}
\end{array}
\right|=\ldots=\left|
\begin{array}{cc}
	\nabla_{N-1}f\textrm{ }\nabla_Nf\\
	q_{ijN-1}\textrm{ }q_{ijN}
\end{array}
\right|=0,
\end{equation}
then exists function $\mu_{ij}$ such that
\begin{equation}
\int_{\partial A} \mu_{ij}df=\sum_{l<m} \int\int_AR_{ijlm}\omega_{l}\wedge\omega_{m}.
\end{equation} 
\\
\textbf{Proof.}\\
Obviously we can write
$$
\frac{\nabla_1f}{q_{ij1}}=\frac{\nabla_2f}{q_{ij2}}=\ldots=\frac{\nabla_Nf}{q_{ijN}}=\frac{1}{\mu_{ij}}.
$$
From Theorem 3 we have
$$
\nabla_l(q_{ijm})-\nabla_m(q_{ijl})+\sum^{N}_{k=1}q_{ijk}Q_{lkm}=R_{ijlm}.
$$
Hence
$$
\nabla_{l}(\mu_{ij}\nabla_mf)-\nabla_{m}(\mu_{ij}\nabla_lf)+\sum^{N}_{k=1}(\nabla_kf)\mu_{ij}Q_{lkm}=R_{ijlm},
$$
or equivalently using Theorem 2
$$
\nabla_mf\cdot\nabla_l \mu_{ij}-\nabla_lf\cdot\nabla_{m}\mu_{ij}=R_{ijlm}.
$$
Hence from relation (61) (Stokes formula) we get the result.\\
\\
\textbf{Note 3.}\\
Condition (62) is equivalent to say that exist functions $\mu_{ij}$ and $f$ such that
\begin{equation}
\omega_{ij}=\mu_{ij}df.
\end{equation}
\\
\textbf{Theorem 10.}\\
If exists function field $F_{ij}$ such that 
\begin{equation}
\sum_{l<m}\sum_{i,j\in I}F_{ij}\Theta^{(1)}_{lm}\left(\overline{x}\omega_{ij}\right)=\overline{0},
\end{equation}
then
\begin{equation}
\overline{x}=\sum^{N}_{k=1}\left(\frac{\sum_{i,j\in I}F_{ij}r_{ijk}}{\sum_{l<m}\sum_{i,j\in I}F_{ij}R_{ijlm}}\right)\overline{e}_k.
\end{equation}
\\
\textbf{Proof.}\\
From (50) and (47) we have
$$
\sum_{l<m}\Theta_{lm}\left(\overline{x}\omega_{ij}\right)=\sum_{l<m}\left|
\begin{array}{cc}
  \nabla_l(\overline{x})\textrm{  }\nabla_m(\overline{x})\\
	q_{ijl}\textrm{  }\textrm{  }q_{ijm}
\end{array}
\right|+\overline{x}\sum_{l<m}R_{ijlm}=
$$
$$
=\sum_{l<m}\left|
\begin{array}{cc}
  \overline{e}_l\textrm{  }\textrm{  }\overline{e}_m\\
	q_{ijl}\textrm{  }\textrm{  }q_{ijm}
\end{array}
\right|+\overline{x}\sum_{l<m}R_{ijlm}=-\sum^{N}_{k,s=1}\epsilon_{ks}q_{ijs}\overline{e}_k+\overline{x}\sum_{l<m}R_{ijlm}=
$$
\begin{equation}
=-\sum^{N}_{k=1}r_{ijk}\overline{e}_k+\overline{x}\sum_{l<m}R_{ijlm},
\end{equation}
where the $r_{ijk}$ are defined in (68) and the $\epsilon_{ks}$ are that of Definition 6 below.
Hence, if exists function field $F_{ij}$ such that 
(65) holds, then we get the validity of (66).\\
\\

Above we have set
\begin{equation}
r_{ijk}=\sum^{N}_{s=1}\epsilon_{ks}q_{ijs}.
\end{equation}
Also
\begin{equation}
\left\langle\sum_{l<m}\Theta_{lm}\left(\overline{x}\omega_{ij}\right),\overline{e}_k\right\rangle=-r_{ijk}+w_k\sum_{l<m}R_{ijlm},
\end{equation}
where $w_k=\left\langle \overline{x},\overline{e}_k\right\rangle$ is called support function of the hypersurface.\\
\\

Setting the symbols 
$$
f_k=f_{k}(x_1,x_2,\ldots,x_N)\textrm{, }k=1,2,\ldots,N
$$
be such that
\begin{equation}
\overline{x}=\sum^{N}_{k=1}f_k\overline{e}_{k},
\end{equation}
then
$$
d(\overline{x})=\sum^{N}_{k=1}df_k\overline{e}_k+\sum^{N}_{k=1}f_kd(\overline{e}_k)=\sum^{N}_{k=1}(df_k)\overline{e}_k+\sum^{N}_{k=1}f_k\sum^{N}_{j=1}\omega_{kj}\overline{e}_j=
$$
$$
=\sum^{N}_{k=1}(df_k)\overline{e}_k+\sum^{N}_{j=1}f_j\sum^{N}_{k=1}\omega_{jk}\overline{e}_k=\sum^{N}_{k=1}\left(df_k+\sum^{N}_{j=1}f_j\omega_{jk}\right)\overline{e}_k.
$$
Hence from (5) and the above relation we get
\begin{equation}
\omega_k=df_k+\sum^{N}_{j=1}f_j\omega_{jk}.
\end{equation}
If we use the Pfaff expansion of the differential we get
$$
\omega_k=\sum^{N}_{l=1}(\nabla_lf_k)\omega_l+\sum^{N}_{j=1}f_j\sum^{N}_{m=1}q_{jkm}\omega_m
$$
Or equivalent
$$
\omega_k=\sum^{N}_{l=1}(\nabla_lf_k)\omega_l+\sum^{N}_{l=1}\sum^{N}_{j=1}f_jq_{jkl}\omega_l
$$
Hence
$$
\nabla_lf_k-\sum^{N}_{j=1}f_jq_{kjl}=\delta_{lk}
$$
Or equivalently we conclude that: 

The necessary conditions such that $f_k$ be the coordinates of the vector $\overline{x}$ (who generates the space $\textbf{S}$), in the moving frame $\overline{e}_{k}$, are 
$$
f_{k;l}=\delta_{kl}.
$$
Hence we get the next\\
\\
\textbf{Theorem 11.}\\
If 
\begin{equation}
f_{k;l}=\delta_{kl}\textrm{, where }k,l\in\{1,2,\ldots,N\},
\end{equation}
then we have
$$
\overline{x}=\sum^{N}_{k=1}f_k\overline{e}_k+\overline{h}\textrm{, where }d\overline{h}=0
$$
and the opposite.\\
\\

Finally having in mind of (66) we get
\begin{equation}
f_k=\frac{\sum_{i,j\in I}F_{ij}r_{ijk}}{\sum_{l<m}\sum_{i,j\in I}F_{ij}R_{ijlm}}
\end{equation}
For two generalized hyper-vectors $F=F_{ij\ldots k}$ and $G=G_{ij\ldots k}$ we define the generalized inner product as
\begin{equation}
(F,G):=\sum^{N}_{i,j,\ldots,k=1}F_{ij\ldots k}G_{ij\ldots k}.
\end{equation}
Hence relation (73) can be written as
$$
f_k=\sum^{N}_{i,j=1}\frac{F_{ij}}{(F,U)}r_{ijk},
$$ 
where 
$$
U:=U_{ij}:=\sum_{l<m}R_{ijlm}
$$
Hence
$$
f_{k;l}=\sum^{N}_{i,j=1}\nabla_l\left(\frac{F_{ij}}{(F,U)}r_{ijk}\right)-\sum^{N}_{i,j,m=1}\frac{F_{ij}}{(F,U)}r_{ijm}q_{kml}=
$$
$$
=\sum^{N}_{i,j=1}\nabla_l\left(\frac{F_{ij}}{(F,U)}\right)r_{ijk}+\sum^{N}_{i,j=1}\frac{F_{ij}}{(F,U)}\nabla_l\left(r_{ijk}\right)-\sum^{N}_{i,j,m=1}\frac{F_{ij}}{(F,U)}r_{ijm}q_{kml}=
$$
$$
=\sum^{N}_{i,j=1}\nabla_l\left(\frac{F_{ij}}{(F,U)}\right)r_{ijk}+\sum^{N}_{i,j=1}\frac{F_{ij}}{(F,U)}r_{ij\{k\};l}
$$
\\
\textbf{Theorem 12.}\\
If $F_{ij}$ is such that
\begin{equation}
\sum^{N}_{i,j=1}\nabla_l\left(\frac{F_{ij}}{(F,U)}\right)r_{ijk}+\sum^{N}_{i,j=1}\frac{F_{ij}}{(F,U)}r_{ij\{k\};l}=\delta_{kl},
\end{equation}
then exists constant vector $\overline{h}$ such that
\begin{equation}
\overline{x}=\sum^{N}_{k,i,j=1}\frac{F_{ij}}{(F,U)}r_{ijk}\overline{e}_k+\overline{h},
\end{equation}
where $d\overline{h}=0$.\\
\\
\textbf{Theorem 13.}\\
There holds (see note 4 below):
\begin{equation}
q_{ik\{m\};l}-q_{ik\{l\};m}=R_{iklm}.
\end{equation}
and
$$
b_{\{k\}m,l}=b_{\{k\}l,m},\eqno{(77.1)}
$$
where 
$$
Y_{\{n\}m,l}=\partial_lY_{nm}-\sum^{N}_{j=1}\Gamma_{njl}Y_{jm}\textrm{, }Y_{\{n\}m,l}=\partial_lY_{nm}-\sum^{N}_{j=1}\Gamma_{mjl}Y_{nj}\eqno{(77.2)}
$$
and
$$
Y_{nm,l}=\partial_lY_{nm}-\sum^{N}_{j=1}\Gamma_{njl}Y_{jm}-\sum^{N}_{j=1}\Gamma_{mjl}Y_{nj}\textrm{, }\ldots\textrm{etc}.\eqno{(77.3)} 
$$
\\
\textbf{Proof.}\\
Easy\\
\\
\textbf{Note 4.}\\
In general in $t_{ij\ldots \{k\}\ldots m;l}$ the brackets indicate where the differential acts. Hence
\begin{equation}
t_{ij\ldots \{k\}\ldots m;l}=\nabla_lt_{ij\ldots k\ldots m}-\sum^{N}_{\nu=1}q_{k\nu l}t_{ij\ldots \nu\ldots m}
\end{equation} 
Also we can use more brackets $\{\}$ in the vector.
$$
t_{ij\ldots \{k_1\}\ldots\{k_2\}\ldots m;l}=
$$
\begin{equation}
=\nabla_lt_{ij\ldots k_1\ldots k_2\ldots  m}-\sum^{N}_{\nu_1=1}q_{k_1\nu_1 l}t_{ij\ldots \nu_1\ldots k_2\ldots m}-\sum^{N}_{\nu_2=1}q_{k_2\nu_2 l}t_{ij\ldots k_1\ldots\nu_2\ldots m}.
\end{equation}
In case we change the ''$;$'' with ''$,$'' then we lead to the classical invariant derivative
$$
t_{ij\ldots \{k\}\ldots m,l}=\partial_lt_{ij\ldots k\ldots m}-\sum^{N}_{\nu=1}\Gamma_{k\nu l}t_{ij\ldots \nu\ldots m},
$$
... etc. The two kinds of derivative lead us to the same evaluations. More precicely it holds
$$
\sum^{N}_{l=1}t_{ij\ldots \{k\}\ldots m;l}\omega_l=\sum^{N}_{l=1}t_{ij\ldots \{k\}\ldots m,l}du_l
$$
\\
\textbf{Theorem 14.}\\
The following forms are invariants of the space $\textbf{S}$:\\
i) The linear element is
\begin{equation}
I=(ds)^2=\sum^{N}_{k=1}\omega_k^2.
\end{equation}
ii) The volume element of $\textbf{S}$ is
\begin{equation}
V=\omega_1\wedge\omega_2\wedge\ldots\wedge\omega_{N-1}\wedge \omega_{N}.
\end{equation}
The area element of the subspace normal to $\overline{e}_{M}$ is
$$
E_{M}=\omega_1\wedge\omega_2\wedge\ldots\wedge\omega_{M-1}\wedge\omega_{M+1}\wedge\ldots\wedge\omega_{N}
$$
iii) The second invariant forms are
\begin{equation}
II_{M}=\left\langle d\overline{x},d\overline{e}_{M}\right\rangle=\sum^{N}_{k=1}\omega_k\omega_{Mk}\textrm{, }M=1,\ldots,N
\end{equation}
iv) The linear element of $\overline{e}_{M}$ is
\begin{equation}
III_{M}=(d\overline{e}_M)^2=\left\langle d\overline{e}_{M},d\overline{e}_{M}\right\rangle=\sum^{N}_{k=1}\omega_{Mk}^2
\end{equation}
v) The Gauss curvature $K_{M}$ which corresponds to the subspace normal to $\overline{e}_M$ vector, is given from 
\begin{equation}
K_{M}=det\left(\kappa^{\{M\}}_{ij}\right)=\frac{\omega_{M1}\wedge\omega_{M2}\wedge\ldots\wedge\omega_{M(M-1)}\wedge\omega_{M(M+1)}\wedge\ldots\wedge \omega_{MN}}{E_{M}},
\end{equation}
where 
\begin{equation}
\kappa^{\{M\}}_{km}:=q_{kMm}.
\end{equation}
The above forms  remain unchanged, in every change of position of $\overline{x}$ and rotation of $\{\overline{e}_j\}_{j=1,2,\ldots,N}$, except for possible change of sign.\\ 
\\
\textbf{Definition 5.}\\
 We define the Beltrami differential operator
\begin{equation}
\Delta_2f:=\sum^{N}_{l<m}\left(
\nabla_{l}^2f+\nabla^2_{m}f+\frac{1}{N-1}\sum^{N}_{i<j}
\left|\begin{array}{cc}
	\nabla_if\textrm{ }\nabla_jf\\
	Q_{lim}\textrm{ }Q_{ljm}
\end{array}\right|\right).
\end{equation} 
We also call $f$ harmonic if $\Delta_2f=0$.\\
\\
\textbf{Remark 2.}\\
In the particular case of $N=2$ where the space $\textbf{S}$ is a two dimensional surface embeded to $\textbf{E}_3$, we have
\begin{equation}
\Delta_2A=\nabla_1\nabla_1A+\nabla_2\nabla_2A+q_2\nabla_1A-q_1\nabla_2A,
\end{equation}
where
\begin{equation}
d\overline{x}=\omega_1\overline{e}_1+\omega_2\overline{e}_2
\end{equation} 
and $d\omega_1=q_1\omega_1\wedge\omega_2$ and $d\omega_2=q_2\omega_1\wedge\omega_2$, there exists functions $f,f^{*}$ such that $\nabla_1f=-\nabla_2f^{*}$ and $\nabla_2f=\nabla_1f^{*}$ 
and $\Delta_2f=\Delta_2f^{*}=0$, ($f,f^{*}$ harmonics). In higher dimensions is not such easy to construct harmonic functions. However we will give here one way of such construction. Before going to study this operator, we simplify the expansion (86) (of the definition of Beltrami differential operator $\Delta_2$). We also generalize it (in a way) as we show below in Definition 8. First we give a definition.\\
\\
\textbf{Definition 6.}\\
Set
\begin{equation}
A_s:=\sum_{l<m}Q_{lsm}\textrm{, }R_k:=\sum^{N}_{s=1}\epsilon_{ks}A_s
\end{equation}
and more generally if
\begin{equation}
a_s=\sum^{N}_{l<m}t_{lsm}\textrm{ and }
r_{s}:=\sum^{N}_{k=1}\epsilon_{sk}a_k,
\end{equation}
where $\epsilon_{ks}:=-1$ if $k<s$, $1$ if $k>s$ and $0$ if $k=s$. Then from identity $\sum^{N_1}_{k,s=1}\epsilon_{ks}f_kf_s=0$, for every choice of $f_k$, we get
\begin{equation}
\sum^{N}_{k=1}R_kA_k=0\textrm{, }\sum^{N}_{s=1}r_s\alpha_s=0
\end{equation}

From the above definition we have the next\\
\\
\textbf{Theorem 15.}
\begin{equation}
\Delta_2f=\sum^{N}_{k=1}\left((N-1)\nabla_k^2f-\frac{1}{N-1}\left
(\nabla_kf\right)R_k\right).
\end{equation}
Also
\begin{equation}
\Delta_2f=(N-1)\sum^{N}_{k=1}\nabla_k^2f-\frac{1}{N-1}\left\langle\overline{\textrm{grad}(f)},\overline{R}\right\rangle.
\end{equation}
\\
\textbf{Proof.}
$$
\Delta_2f=\sum_{l<m}\left(\nabla_l^2f+\nabla_m^2f-\frac{1}{N-1}\sum^{N}_{k,s=1}\left(\nabla_kf\right)\epsilon_{ks}Q_{lsm}\right)=
$$
$$
=(N-1)\sum^{N}_{k=1}\nabla_k^2f-\frac{1}{N-1}\sum^{N}_{k=1}\left(\nabla_kf\right)\sum^{N}_{s=1}\epsilon_{ks}A_s=
$$
$$
=(N-1)\sum^{N}_{k=1}\nabla_k^2f-\frac{1}{N-1}\sum^{N}_{k=1}\left(\nabla_kf\right)R_k.
$$
\\
\textbf{Note 5.}\\
1) If exists $f$ such $\sum^{N}_{k=1}A_k\omega_k=df$, then $\nabla_kf=A_k$ and
$$
\Delta_2f=(N-1)\sum_{k=1}^N\nabla_kA_k.
$$
2) If $h_k$ is any vector such that $\sum^{N}_{k=1}h_k=0$
and $f_k=\nabla_kF$ is such that
\begin{equation}
(N-1)\nabla_kf_k-\frac{1}{N-1}f_kR_k=h_k\textrm{, }\forall k=1,2,\ldots,N,
\end{equation}
then
$$
\Delta_2F=0.
$$
Hence solving first the PDE's
$$
(N-1)\nabla_k\Psi-\frac{1}{N-1}\Psi R_k=h_k\textrm{, }\forall k=1,2,\ldots,N,\eqno{(eq1)}
$$
$$
\sum^{N}_{k=1}h_k=0,\eqno{(eq2)}
$$
we find $N$ solutions $\Psi=f_k$. Then we solve $f_k=\nabla_kF$ and thus we get a solution $F$ of $\Delta_2F=0$.\\
3) This lead us to define the derivative $D_k$, $k=1,2,\ldots,N$, which acts in every function $f$ as
\begin{equation}
D_k\left(f\right):=D_kf:=(N-1)\nabla_kf-\frac{1}{N-1}R_k f.
\end{equation}
Obviously $D_k$ is linear. But
$$
d(fg)=gdf+fdg=g\sum^{N}_{k=1}\nabla_kf\omega_k+f\sum^{N}_{k=1}\nabla_kg\omega_k=\sum^{N}_{k=1}\left(g\nabla_kf+f\nabla_kg\right)\omega_k,
$$
or equivalently
\begin{equation}
\nabla_k(fg)=g\nabla_kf+f\nabla_kg.
\end{equation}
Hence the differential operator $D_k$ acts in the product of $f$ and $g$ as
$$
D_k(fg)=(N-1)\nabla_k(fg)-\frac{1}{N-1}fgR_k=
$$
$$
=(N-1)f\nabla_kg+(N-1)g\nabla_kf-\frac{1}{N-1}fgR_k=fD_kg+(N-1)g\nabla_kf=
$$  
$$
=gD_kf+(N-1)f\nabla_kg
$$
and finally we have
\begin{equation}
D_k(fg)=fD_kg+gD_kf+\frac{1}{N-1}f gR_k,
\end{equation}
\begin{equation}
fD_kg-gD_kf=(N-1)(f\nabla_kg-g\nabla_kf).
\end{equation}
The commutator of $D_k$ and $\nabla_k$ acting to a scalar field is
$$
\left[D_k,\nabla_k\right]f
=D_k\nabla_kf-\nabla_kD_kf=
$$
$$
=(N-1)\nabla_k^2f-\frac{1}{N-1}R_k\nabla_kf-\left((N-1)\nabla_k^2f-\frac{1}{N-1}\nabla_k(R_kf)\right).
$$
Hence finally after simplifications
\begin{equation}
\left[D_k,\nabla_k\right]f
=\frac{\nabla_kR_k}{N-1}f.
\end{equation}
i) If $f=const$, then 
$$
D_kf=-\frac{1}{N-1}R_kf.
$$
ii) The derivative $D_k$ is such that if exists function $f$ with 
$$
D_kf=0\textrm{, }\forall k=1,2,\ldots,N,
$$ 
then equivalently
$$
\nabla_k(\log f)=\frac{1}{(N-1)^2}R_k\textrm{, }\forall k=1,2,\ldots,N.
$$
But then
\begin{equation}
D_kf=0\Leftrightarrow \exists g\left(=(N-1)^2\log f\right):R_k=\nabla_kg\Leftrightarrow R=dg.
\end{equation}
By these arguments we conclude to the result, that only in specific spaces $\textbf{S}$ exists $f$ such $D_kf=0$, for all $k=1,2,\ldots,N$. Hence we have a definition-theorem:\\
\\
\textbf{Definition 7.}\\
A space is called $\textbf{S}$-$R$ iff exists function $g$ such that $R=dg$.\\
\\
\textbf{Theorem 16.}\\
A space $\textbf{S}$ is $\textbf{S}$-$R$ iff exists $g$ such that $D_kg=0$, for all $k=1,2,\ldots,N$.\\
\\
\textbf{Corollary 2.}\\
In a $\textbf{S}$-$R$ space
\begin{equation}
\Theta_{lm}(R)=0\textrm{, }\forall l,m\in\{1,2,\ldots,N\}.
\end{equation}
\textbf{Proof.}\\
It is $R=\sum^{N}_{k=1}R_k\omega_k$ and $\textbf{S}$ is $\textbf{S}$-$R$. Hence
$$
\Theta_{lm}(R)=\nabla_lR_m-\nabla_mR_l+\sum^{N}_{k=1}R_kQ_{lkm}.
$$    
But exists $g$ such that $R_k=\nabla_kg$, for all $k=1,2,\ldots,N$. Hence
$$
\Theta_{lm}(R)=\nabla_l\nabla_mg-\nabla_m\nabla_lg+\sum^{N}_{k=1}(\nabla_kg)Q_{lkm}=0.
$$
The last equality is due to Theorem 2. Hence the result follows.\\     
\\
\textbf{Corollary 3.}\\
In a $\textbf{S}$-$R$ space holds 
$$
d(fR)=\sum_{l<m}\left(R_m\nabla_lf-R_l\nabla_mf \right)\omega_l\wedge\omega_m.\eqno{(101.1)}
$$
\textbf{Proof.}\\
The result is application of Theorem 7.\\
\\
\textbf{Corollary 4.}\\
In every space $\textbf{S}$ holds
$$
\sum^{N}_{k=1}A_k(D_kf)=(N-1)\sum^{N}_{k=1}A_k(\nabla_kf).\eqno{(101.2)}
$$
\\
\textbf{Theorem 17.}\\
Set  
$$
Df:=(D_1f)\omega_1+(D_2f)\omega_2+\ldots+(D_Nf)\omega_N,\eqno{(101.3)}
$$
then
$$
Df=(N-1)df-\frac{1}{N-1}fR.\eqno{(101.4)}
$$
Also
$$
D\left\langle\overline{V}_1,\overline{V}_2\right\rangle=\left\langle D\overline{V}_1,\overline{V}_2\right\rangle+\left\langle \overline{V}_1,D\overline{V}_2\right\rangle+\frac{R}{N-1}\left\langle \overline{V}_1,\overline{V}_2\right\rangle.\eqno{(101.5)}
$$
\\
\textbf{Theorem 18.}\\
If a space $\textbf{S}$ is $\textbf{S}$-$R$, then the PDE $\Delta_2\Psi=\mu \Psi$ have non trivial solution.\\
\\
\textbf{Proof.}\\
If $\textbf{S}$ is $\textbf{S}$-$R$, then from (100) we get that exists $f,g$ such $g=(N-1)^2\log f$ and $R=dg$,  $D_jf=0$, $\forall j$. Hence 
$$
(N-1)\nabla_jf-\frac{1}{N-1}R_jf=0\Rightarrow
$$
$$
(N-1)\nabla_j\nabla_jf-\frac{1}{N-1}R_j\nabla_jf-\frac{\nabla_jR_j}{N-1}f=0\Rightarrow
$$
$$
\Delta_2f=\frac{f}{N-1}\sum^{N}_{j=1}\nabla_jR_j.
$$
\\
\textbf{Theorem 19.}\\
In every space $\textbf{S}$ we have
\begin{equation}
\overline{Df}:=\sum^{N}_{i=1}(D_kf)\overline{e}_k=(N-1)\overline{\textrm{grad}f}-\frac{1}{N-1}\overline{R}f.
\end{equation}
\\
\textbf{Note 6.}\\
1) If $C$ is a curve, then
$$
\left(\frac{Df}{ds}\right)_C=(N-1)\left(\frac{df}{ds}\right)_C-\frac{1}{N-1}\left(\frac{dR}{ds}\right)_C f
$$
and if $(dR)_C\neq0$, then
$$
\left(\frac{Df}{dR}\right)_C=(N-1)\left(\frac{df}{dR}\right)_C-\frac{1}{N-1}f
$$
and $\left(\frac{Df}{dR}\right)_C=0$ iff 
$$
\left(\frac{df}{dR}\right)_C=\frac{1}{(N-1)^2}f.
$$
Hence
$$
f_C=const\cdot\exp\left(\frac{R_C}{(N-1)^2}\right).
$$
Hence if $f=f(x_1,x_2,\ldots,x_N)$ and $C:x_i=x_i(s)$, $s\in\left[a,b\right]$, such $\left(\frac{Df}{ds}\right)_C=0$, then 
$$
f\left(x_1(s),x_2(s),\ldots,x_N(s)\right)=const\cdot\exp\left(\frac{R\left(x_1(s),x_2(s),\ldots,x_N(s)\right)}{(N-1)^2}\right).
$$
2) Also assuming $C:x_i=x_i(s)$ is a curve in $\textbf{S}$ and $\overline{\nu}$ is the tangent vector of $C$, which is such that $\left\langle \overline{\nu},\overline{\nu}\right \rangle=1\Rightarrow\left\langle \overline{\nu},\frac{d\overline{\nu}}{ds}\right \rangle=0$, then
$$
\frac{D\overline{\nu}}{ds}=(N-1)\frac{d\overline{\nu}}{ds}-\frac{1}{N-1}\frac{dR}{ds}\overline{\nu}.
$$ 
Consequently we get
$$
\left\langle\frac{D\overline{\nu}}{ds},\frac{d\overline{\nu}}{ds}\right\rangle=(N-1)k(s)
$$
and
$$
\left\langle\frac{D\overline{\nu}}{ds},\overline{\nu}\right\rangle=-\frac{1}{N-1}\frac{dR}{ds}.
$$
But
$$
D\left\langle \overline{\nu},\overline{\nu}\right \rangle=D1\Leftrightarrow 2\left\langle D\overline{\nu},\overline{\nu}\right \rangle+\frac{dR}{N-1}=-\frac{dR}{N-1}
.
$$
Hence
\begin{equation}
\left\langle D\overline{\nu},\overline{\nu}\right \rangle=-\frac{dR}{N-1}.
\end{equation}
Hence differentiating we get
$$
d\left\langle D\overline{\nu},\overline{\nu}\right \rangle=0.
$$
\\
\textbf{Theorem 20.}\\
For every curve $C:x_i=x_i(s)$ of a general space $\textbf{S}$, with $s$ being its normal parameter, the unitary tangent vector $\overline{\nu}$ of $C$ has the property
$$
\left\langle \left(\frac{D\overline{\nu}}{ds}\right)_C,\overline{\nu}\right\rangle=-\frac{1}{N-1}\left(\frac{dR}{ds}\right)_C.
$$
Moreover we define 
$$
k^{*}=\left(k^{*}\right)_C:=\left(\frac{1}{\rho_{R}}\right)_C:=\left|\left(\frac{D \overline{\nu}}{ds}\right)_C\right|,
$$
then
$$
k^{*}(s)=\sqrt{(N-1)^2k^2(s)+\frac{1}{(N-1)^2}\left(\left(\frac{dR}{ds}\right)_C\right)^2}.
$$
If the space is $\textbf{S}$-$R$, then
$$
\left(k^{*}\right)_C=(N-1)k(s).
$$
\\
\textbf{Corollary 5.}\\
Let $C:x_i=x_i(s)$ be a curve of a $\textbf{S}$-$R$ space. If $s$ is the normal parameter of $C$ and $\overline{\nu}$ is the unitary tangent vector of $C$, then
$$
\frac{D\overline{\nu}}{ds}=(N-1)\frac{d\overline{\nu}}{ds}\textrm{ 
 and
 }\left|\frac{D\overline{\nu}}{ds}\right|=(N-1)k(s),
$$
where $k(s)$ is the curvature of $C$.\\
\\

Next we generalize the deffinition of $\Delta_2$ operator:\\
\\
\textbf{Definition 8.}\\
\begin{equation}
\Delta^{(\lambda)}(f)=\sum^{N}_{i,j=1}\lambda_{ij}D_i\nabla_jf.
\end{equation}
Hence if $\lambda_{ij}=\delta_{ij}$, then 
\begin{equation}
\Delta_2f=\sum^{N}_{i,j=1}\delta_{ij}D_i\nabla_jf=\sum^{N}_{k=1}D_k\nabla_kf.
\end{equation}
\\
\textbf{Theorem 21.}\\
If $\lambda_{ij}=\epsilon_{ij}$, then
\begin{equation}
\Delta^{(\epsilon)}f=\sum^{N}_{k=1}\left((N-1)A_k+\frac{1}{N-1}\sum^{N}_{s=1}\epsilon_{ks}R_s\right)\nabla_kf.
\end{equation}
\\
\textbf{Proof.}\\
For $\lambda_{ij}=\epsilon^{*}_{ij}=-\epsilon_{ij}$ we have
$$
\Delta^{(\epsilon^{*})}f=\sum^{N}_{i,j=1}\epsilon^{*}_{ij}D_i\nabla_jf=\sum_{i<j}\epsilon^{*}_{ij}D_i\nabla_{j}f+\sum_{i>j}\epsilon^{*}_{ij}D_i\nabla_jf=
$$
$$
=\sum_{i<j}\left(D_i\nabla_jf-D_j\nabla_if\right)=
$$
$$
=(N-1)\sum_{i<j}\left(\nabla_i\nabla_jf-\nabla_j\nabla_if\right)-\frac{1}{N-1}\sum_{i<j}\left(R_i\nabla_jf-R_j\nabla_if\right)=
$$
$$
=-(N-1)\sum_{i<j}\sum^{N}_{k=1}(\nabla_kf)Q_{ikj}-\frac{1}{N-1}\sum^{N}_{k,s=1}\epsilon_{ks}R_s\nabla_kf=
$$
$$
=-(N-1)\sum^{N}_{k=1}A_k\nabla_kf-\frac{1}{N-1}\sum^{N}_{k,s=1}\epsilon_{ks}R_s\nabla_kf=
$$
$$
=-\sum^{N}_{k=1}\left((N-1)A_k+\frac{1}{N-1}\sum^{N}_{s=1}\epsilon_{ks}R_s\right)\nabla_kf.
$$
$qed$\\ 
\\
\textbf{Example 1.}\\
If $\lambda_{ij}=\epsilon_{ij}$ and $\textbf{S}$ is $\textbf{S}$-$R$, then 
\begin{equation}
\Delta^{(\epsilon)}g=0,
\end{equation}
where $dg=R$. That is because
$$
\Delta^{(\epsilon)}g=(N-1)\sum^{N}_{k=1}A_kR_k+\frac{1}{N-1}\sum^{N}_{k,s=1}\epsilon_{ks}R_kR_s=0,
$$
since $\sum^{N}_{k=1}A_kR_k=0$ and $\sum^{N}_{k,s=1}\epsilon_{ks}R_kR_s=0$. Also easily we get
\begin{equation}
\left\langle\Delta^{(\epsilon)}\overline{x},\overline{R}\right\rangle=0.
\end{equation}
\\
\textbf{Example 2.}\\
In a $\textbf{S}$-$R$ space we have $R=dg$, for a certain $g$ and thus $\nabla_kg=R_k$. Hence we can write 
\begin{equation}
\Delta_2g=(N-1)\sum^{N}_{k=1}\nabla_kR_k-\frac{1}{N-1}\sum^{N}_{k=1}R_k^2.
\end{equation}
\\
\textbf{Theorem 22.}\\
In the general case when $\lambda_{ij}$ is any field we 
can write
\begin{equation}
\Delta^{(\lambda)}f=(N-1)\sum^{N}_{i,j=1}\lambda_{ij}\nabla_{i}\nabla_{j}f-\frac{1}{N-1}\sum^{N}_{i,j=1}\lambda_{ij}R_{i}\nabla_{j}f.
\end{equation}
and
\begin{equation}
\left\langle\Delta^{(\lambda)}\overline{x},\overline{e}_k\right\rangle=(N-1)\sum_{i,j=1}^N\lambda_{ij}q_{jki}-\frac{1}{N-1}\sum^{N}_{i=1}\lambda_{ik}R_i.
\end{equation}
\\
\textbf{Proof.}\\
The first identity is easy. For the second we have:
Setting $f\rightarrow \overline{x}$, we get
$$
\Delta^{(\lambda)}\overline{x}=\sum^{N}_{i,j=1}\lambda_{ij}D_i\left(\nabla_j\overline{x}\right)=\sum^{N}_{i,j=1}\lambda_{ij}D_i\left(\overline{e}_j\right)=
$$
$$
=\sum^{N}_{i,j=1}(N-1)\lambda_{ij}\nabla_i\overline{e}_j-\sum^{N}_{i,j=1}\frac{\lambda_{ij}}{N-1}R_i\overline{e}_j=
$$
$$
=\sum^{N}_{i,j=1}\lambda_{ij}\left((N-1)\nabla_i\overline{e}_j-\frac{1}{N-1}R_i\overline{e}_j\right)=
$$
$$
=\sum^{N}_{i,j,k=1}(N-1)\lambda_{ij}q_{jki}\overline{e}_k-\sum^{N}_{i,k=1}\frac{\lambda_{ik}}{N-1}R_i\overline{e}_k
$$
$$
=\sum^{N}_{k=1}\left((N-1)\sum^{N}_{i,j=1}\lambda_{ij}q_{jki}-\frac{1}{N-1}\sum^{N}_{i=1}\lambda_{ik}R_i\right)\overline{e}_k.
$$
\\

We need some notation to proceed further\\
\\
\textbf{Definition 9.}\\
For any field $\lambda_{ij}$ we define 
$$
R^{\{\lambda\}}_k:=\sum^{N}_{i=1}\lambda_{ik}R_{i}\eqno{(111.1)}
$$
$$
A^{(\lambda)}_{k}:=\sum_{i<j}\lambda_{ij}Q_{ikj}\textrm{, }A^{{(\lambda)}{*}}_k:=\sum_{i<j}\lambda_{ij}Q^{*}_{ikj},\eqno{(111.2)}
$$
where
$$
Q_{ikj}:=q_{ikj}-q_{jki}\textrm{, }Q^{*}_{ikj}:=q_{ikj}+q_{jki}\eqno{(111.3)}
$$
and
$$
R^{(\lambda)}_{k}:=\sum^{N}_{s=1}\epsilon_{ks}A^{(\lambda)}_s\textrm{, }R^{{(\lambda)}{*}}_{k}:=\sum^{N}_{s=1}\epsilon_{ks}A^{{(\lambda)}{*}}_s.\eqno{(111.4)}
$$
\\
\textbf{Proposition 5.}
$$
\sum^{N}_{s=1}A^{(\lambda)}_sR^{(\lambda)}_s=0\textrm{, }\sum^{N}_{s=1}A^{{(\lambda)}{*}}_sR^{{(\lambda)}{*}}_s=0.\eqno{(111.5)}
$$
\\
\textbf{Theorem 23.}\\
Let $\lambda_{ij}$ be antisymmetric i.e. $\lambda_{ij}=-\lambda_{ji}$, then we have in general:\\ 
1)
$$
\Delta^{(\lambda)}f=-\sum^{N}_{k=1}\left((N-1)A^{(\lambda)}_k+\frac{1}{N-1}\sum^{N}_{i=1}\lambda_{ik}R_i\right)\nabla_kf.
$$
2)
$$
\left\langle\Delta^{(\lambda)}\overline{x},\overline{e}_k\right\rangle=-(N-1)A^{(\lambda)}_k-\frac{1}{N-1}\sum^{N}_{i=1}\lambda_{ik}R_i.
$$
3) 
$$
\sum^{N}_{k=1}R_kR^{\{\lambda\}}_k=0\Leftrightarrow \left\langle \overline{R},\overline{R}^{\{\lambda\}}\right\rangle=0.
$$
Hence
$$
\left\langle \Delta^{(\lambda)}\overline{x},\overline{R}\right\rangle=-(N-1)\sum^{N}_{k=1}A^{(\lambda)}_kR_k
$$
and
$$
\left\langle \Delta^{(\lambda)}\overline{x},\overline{R}^{(\lambda)}\right\rangle=-\frac{1}{N-1}\sum^{N}_{i,k=1}\lambda_{ik}R_iR^{(\lambda)}_k.
$$
4) In case the space is $\textbf{S}$-$R$, with $\nabla_kg=R_k$, then
$$
\Delta^{(\lambda)}g=-(N-1)\sum^{N}_{k=1}A^{(\lambda)}_kR_k=\left\langle \Delta^{(\lambda)}\overline{x},\overline{R}\right\rangle.
$$
\\
\textbf{Proof.}\\
The case of (1) follows with straight forward evaluation. Since $\lambda_{ij}$ is antisymmetric, we have
$$
\Delta^{(\lambda)}f=\sum_{i<j}\lambda_{ij}D_i\nabla_jf-\sum_{i<j}\lambda_{ij}D_j\nabla_if=
$$
$$
=\sum_{i<j}\lambda_{ij}\left[(N-1)\nabla_i\nabla_jf-\frac{1}{N-1}R_i\nabla_jf-\left((N-1)\nabla_j\nabla_if-\frac{1}{N-1}R_j\nabla_if\right)\right]
$$
$$
=(N-1)\sum_{i<j}\lambda_{ij}\left(\nabla_i\nabla_jf-\nabla_j\nabla_if\right)-\frac{1}{N-1}\sum_{i<j}\lambda_{ij}\left(R_i\nabla_jf-R_j\nabla_if\right)=
$$
$$
=(N-1)\sum_{i<j}\lambda_{ij}\left(-\sum^{N}_{k=1}\nabla_kfQ_{ikj}\right)-\frac{1}{N-1}\sum_{i<j}\lambda_{ij}rot_{ij}(R_i\nabla_jf)=
$$
$$
=-(N-1)\sum^{N}_{k=1}A^{(\lambda)}_k\nabla_kf-\frac{1}{N-1}\sum_{i<j}\lambda_{ij}rot_{ij}\left(R_i\nabla_jf\right)=
$$
$$
=-(N-1)\sum^{N}_{k=1}A^{(\lambda)}_k\nabla_kf-\frac{1}{N-1}\sum^{N}_{k=1}\left(\sum^{N}_{i=1}\lambda_{ik}R_i\right)\nabla_kf=
$$
$$
=-\sum^{N}_{k=1}\left((N-1)A^{(\lambda)}_k+\frac{1}{N-1}\sum^{N}_{i=1}\lambda_{ik}R_i\right)\nabla_kf.
$$
In the second we use $\nabla_k\overline{x}=\overline{e}_k$.\\
The third relation can be proved if we consider the formula
$$
\sum^{N}_{i,k=1}\lambda_{ik}R_iR_k=\sum_{i<k}\lambda_{ik}R_{i}R_k-\sum_{i<k}\lambda_{ik}R_kR_i=0.
$$
The fourth case follows easily from the first with $\nabla_kg=R_k$. Note that in the most general case where $\lambda_{ij}$ is arbitrary we still have
$$
\sum^{N}_{k=1}A^{(\lambda)}_kR^{(\lambda)}_k=0.
$$
\\
\textbf{Theorem 24.}\\
In case $\lambda_{ij}$ is any antisymmetric field, we have
\begin{equation}
\Delta^{(\lambda)}f=\sum^{N}_{k=1}\left\langle\Delta^{(\lambda)}\overline{x},\overline{e}_{k}\right\rangle\nabla_kf=\left\langle\Delta^{(\lambda)}\overline{x},\overline{\textrm{grad}(f)}\right\rangle,
\end{equation}
where
\begin{equation}
\overline{\textrm{grad}(f)}:=\sum^{N}_{k=1}(\nabla_kf)\overline{e}_k.
\end{equation}
\\
\textbf{Proof.}\\
If $\lambda_{ij}$ is any antisymmetric field, then using Theorems 22, 23 (see also and Theorem 2 and Corollary 1) we get the result.\\
\\
\textbf{Corollary 6.}\\
If $\lambda_{ij}$ is any antisymmetric field such that 
\begin{equation}
\Delta^{(\lambda)}\left(\overline{x}\right)=\overline{0},
\end{equation}
then for every $f$ we have
\begin{equation}
\Delta^{(\lambda)}f=0.
\end{equation}
\\
\textbf{Remark 3.}\\The above corollary is very strange. The Beltrami operator over an antisymmetric field is zero for every function $f$ when (114) holds. This force us to conclude that in any space the vector
\begin{equation}
\overline{\sigma}^{(\lambda)}:=\sum^{N}_{k=1}\sigma^{(\lambda)}_k\overline{e}_k=\Delta^{(\lambda)}\left(\overline{x}\right)
\end{equation} 
must play a very prominent role in the geometry of $\textbf{S}$. Then (in any space):
\begin{equation}
\sigma^{(\lambda)}_k=(N-1)\sum_{i,j=1}^N\lambda_{ij}q_{jki}-\frac{1}{N-1}\sum^{N}_{i=1}\lambda_{ik}R_i.
\end{equation}
In case $\lambda_{ij}$ is antisymmetric, then we get
$$
\sigma^{(\lambda)}_k=-(N-1)A_k^{(\lambda)}-\frac{1}{N-1}\sum^{N}_{i=1}\lambda_{ik}R_i=
$$
$$
=-(N-1)A_k^{(\lambda)}-\frac{1}{N-1}\sum^{N}_{l=1}\epsilon^{(\lambda)}_{kl}A_l,
$$
where
\begin{equation}
\epsilon^{(\lambda)}_{kl}:=\sum^{N}_{i=1}\lambda_{ik}\epsilon_{il}.
\end{equation}
Hence when $\lambda_{ij}$ is antisymmetric, then
$$
\sigma^{(\lambda)}_{k}=-(N-1)A_k^{(\lambda)}-\frac{1}{N-1}\sum^{N}_{l=1}\epsilon^{(\lambda)}_{kl}A_l.
$$
Hence
\begin{equation}
\sigma^{(\lambda)}_{k}=-\sum^{N}_{i<j}\left[(N-1)\lambda_{ij}Q_{ikj}+\frac{1}{N-1}\sum^{N}_{l=1}\epsilon^{(\lambda)}_{kl}Q_{ilj}\right].
\end{equation}
\\
\textbf{Definition 10.}\\
We define
\begin{equation}
\Lambda^{(\lambda)}_{kij}:=(N-1)\lambda_{ij}Q_{ikj}+\frac{1}{N-1}\sum^{N}_{l=1}\epsilon^{(\lambda)}_{kl}Q_{ilj}
\end{equation}
and
\begin{equation}
M^{(\lambda)}_{kij}:=(N-1)\lambda_{ij}Q^{*}_{ikj}-\frac{1}{N-1}\sum^{N}_{l=1}\epsilon^{(\lambda)}_{kl}Q_{ilj}.
\end{equation}
Also we define the mean curvature with respect to the $\lambda_{ij}$ as
\begin{equation}
H^{(\lambda)}_k:=(N-1)\sum^{N}_{i,j=1}\lambda_{ij}q_{jki}.
\end{equation}
Then we have
\begin{equation}
H^{(\lambda^{+})}_k=\frac{N-1}{2}\sum^{N}_{i,j=1}\lambda_{ij}Q^{*}_{ikj}
\end{equation}
and
\begin{equation}
H^{(\lambda^{-})}_k=-\frac{N-1}{2}\sum^{N}_{i,j=1}\lambda_{ij}Q_{ikj},
\end{equation}
where $\lambda^{+}$ is the symmetric part of $\lambda_{ij}$ and $\lambda^{-}$ is the antisymmetric part of $\lambda_{ij}$.\\
\\
\textbf{Remark.}\\
It holds
\begin{equation}
\Lambda_{kij}^{(\lambda)}+M_{kij}^{(\lambda)}=2(N-1)\lambda_{ij}q_{ikj}
\end{equation}
\\
\textbf{Theorem 24.1}\\
In any space and any $\lambda_{ij}$, we have
$$
\Delta^{(\lambda)}\left(\overline{x}\right)=\sum^{N}_{k=1}\left(H^{(\lambda)}_k-\frac{1}{N-1}R^{\{\lambda\}}_k\right)\overline{e}_k.\eqno{(125.1)}
$$
Hence also
$$
\Delta^{(\lambda)}\left(\overline{x}\right)=\overline{H}^{(\lambda)}-\frac{1}{N-1}\overline{R}^{\{\lambda\}},\eqno{(125.2)}
$$
where (from Deffinition 9):
$$
\overline{R}^{\{\lambda\}}=\sum^{N}_{k=1}R^{\{\lambda\}}_k\overline{e}_k.
$$
\\
\textbf{Proof.}\\
Easy from the above.\\
\\
\textbf{Theorem 25.}\\
If $\lambda_{ij}$ is antisymmetric we have
\begin{equation}
\sigma^{(\lambda)}_{k}=-\sum^{N}_{i<j}\Lambda^{(\lambda)}_{kij}.
\end{equation}
If $\lambda_{ij}=\lambda^{+}_{ij}+\lambda^{-}_{ij}$, then
$$
\sigma^{(\lambda)}_{k}=\sum^{*}_{i\leq j}M^{(\lambda^+)}_{kij}-\sum_{i<j}\Lambda^{(\lambda^-)}_{kij}=
$$
\begin{equation}
=\sum_{i<j}M^{(\lambda^+)}_{kij}-\sum_{i<j}\Lambda^{(\lambda^-)}_{kij}+(N-1)\sum^{N}_{i=1}\lambda_{ii}q_{iki},
\end{equation}
where the asterisk on the summation means that when $i=j$ we must multiply the summands with $\frac{1}{2}$.\\
\\
\textbf{Corollary 7.}\\
1) We have
\begin{equation}
\Delta^{(\lambda)}\overline{x}=\overline{0},
\end{equation}
iff for all $k=1,2,\ldots,N$ we have
\begin{equation}
H^{(\lambda)}_k-\frac{1}{N-1}R^{\{\lambda\}}_k=0.
\end{equation}
2) For every $\lambda_{ij}$ we have
$$
R^{\{\lambda\}}_k=\sum^{N}_{l,s=1}\epsilon_{ls}\lambda_{lk}\sum_{i<j}Q_{isj}.\eqno{(129.1)}
$$
\\
\textbf{Proof.}\\Actually then we have
\begin{equation}
\sigma_k^{(\lambda)}=H^{(\lambda)}_k-\frac{1}{N-1}\sum^{N}_{i=1}\lambda_{ik}R_i,
\end{equation}
where
\begin{equation}
H^{(\lambda)}_k=(N-1)\sum^{N}_{i,j=1}\lambda_{ij}q_{jki}=H^{(\lambda^{+})}_k+H^{(\lambda^{-})}_k.
\end{equation}
In case $\lambda_{ij}=g_{ij}$ is the metric tensor, then we write
\begin{equation}
H_k:=H^{(g)}_k=(N-1)\sum^{N}_{i,j=1}g_{ij}q_{ikj}
\end{equation}
and call $H_k$ mean curvature of the surface $\textbf{S}$.\\
\\
\textbf{Theorem 26.}\\
If $\lambda_{ij}=g_{ij}$ is the metric tensor, then $\Delta^{(g)}\overline{x}=\overline{0}$ iff
\begin{equation}
H_{k}-\frac{1}{N-1}\sum^{N}_{i=1}g_{ik}R_i=0.
\end{equation}
\\
\textbf{Corollary 8.}\\
If $\lambda_{ij}$ is antisymmetric, then
$$
\Delta^{(\lambda)}(fg)=f\Delta^{(\lambda)}g+g\Delta^{(\lambda)}f
$$
\\
\textbf{Proof.}\\
Use Theorem 24 with
$$
\overline{\textrm{grad}(fg)}=f\overline{\textrm{grad}(g)}+g\overline{\textrm{grad}(f)}.
$$
\\
\textbf{Theorem 27.}\\
Assume that $\lambda_{ij}=\epsilon_{ij}-$antisymmetric. Then 
\begin{equation}
\Delta^{(\epsilon)}\overline{x}=\overline{0}\Leftrightarrow (N-1)A_k+\frac{1}{N-1}\sum^{N}_{l=1}\epsilon^{(2)}_{kl}A_l=0,
\end{equation}
where 
\begin{equation}
\epsilon^{(2)}_{kl}=\sum^{N}_{i=1}\epsilon_{ik}\epsilon_{il}.
\end{equation}
For every space $\textbf{S}$ with $A_{i}$ as above and for  every function $f$, we have
\begin{equation}
\Delta^{(\epsilon)}f=0.
\end{equation}
\\
\textbf{Proof.}\\
Use Theorem 23 with $\lambda_{ij}=\epsilon_{ij}$ and then Definitions 9,6.\\   
\\
\textbf{Remark 5.} From the above propositions we conclude that in every space $\textbf{S}$ we have at least one antisymmetric field (the $\lambda_{ij}=\epsilon_{ij}$) that under condition
\begin{equation} (N-1)A_k+\frac{1}{N-1}\sum^{N}_{l=1}\epsilon^{(2)}_{kl}A_l=0\Leftrightarrow \sum^{N}_{k=1}\Lambda^{(\epsilon)}_{kij}=0,
\end{equation}
we have
\begin{equation}
\Delta^{(\epsilon)}\left(f\right)=0\textrm{, }\forall f. 
\end{equation}
Hence in every space $\textbf{S}$ the quantity
\begin{equation}
\sigma^{(\epsilon)}_k:=(N-1)A_k+\frac{1}{N-1}\sum^{N}_{l=1}\epsilon^{(2)}_{kl}A_l\textrm{, }k=1,2,\ldots,N
\end{equation}
is important. More general the quantities $\sigma^{(\lambda)}_k$ of (119) with $\lambda_{ij}$ antisymmetric are of extreme interest.\\
\\  
\textbf{Corollary 9.}\\
If $\lambda_{ij}$ is antisymmetric, then
\begin{equation}
\Delta^{(\lambda)}f=\frac{1}{(N-1)^2}\left\langle \Delta^{(\lambda)}\overline{x},\overline{R}\right\rangle f+\frac{1}{N-1}\left\langle \overline{Df},\Delta^{(\lambda)}\overline{x}\right\rangle.
\end{equation}
\\
\textbf{Proof.}\\
It holds
$$
\overline{Df}=(N-1)\overline{\textrm{grad} f}-\frac{1}{N-1}\overline{R}f.
$$
Hence
$$
\left\langle\overline{Df},\Delta^{(\lambda)}\overline{x}\right\rangle=(N-1)\left\langle\overline{\textrm{grad}f},\Delta^{(\lambda)}\overline{x}\right\rangle-\frac{1}{N-1}\left\langle \Delta^{(\lambda)}\overline{x},\overline{R}\right\rangle f.
$$
Now since $\lambda_{ij}$ is antisymmetric we have from Theorem (24) the result.\\ 
\\
\textbf{Corollary 10.}\\
If $\lambda_{ij}$ is antisymmetric, then 
\begin{equation}
\left\langle\overline{Df},\Delta^{(\lambda)}\overline{x}\right\rangle=0\Leftrightarrow\Delta^{(\lambda)}f=\frac{1}{(N-1)^2}\left\langle \Delta^{(\lambda)}\overline{x},\overline{R}\right\rangle f.
\end{equation}
In particular
\begin{equation}
\overline{Df}=\overline{0}\Rightarrow\Delta^{(\lambda)}f=\frac{1}{(N-1)^2}\left\langle \Delta^{(\lambda)}\overline{x},\overline{R}\right\rangle f.
\end{equation}
\\
\textbf{Corollary 11.}\\If in a space $\textbf{S}$, the $\lambda_{ij}$ is antisymmetric with
\begin{equation}
\left\langle \Delta^{(\lambda)}\overline{x},\overline{R}\right\rangle=0,
\end{equation}
then 
$$ 
\left\langle \overline{H}^{(\lambda)}, \overline{R}\right\rangle=0
$$
and for every $f$ holds 
\begin{equation}
\overline{Df}=\overline{0}\Rightarrow \Delta^{(\lambda)}f=0.
\end{equation}
\\
\textbf{Remark.}\\
If $\textbf{S}$ is $S$-$R$, then from Example 1, pg.23, we have  $\left\langle \Delta^{(\epsilon)}\overline{x},\overline{R}\right\rangle=0$. Hence in a $S$-$R$ space we have
$$
\Delta^{(\epsilon)}f=\frac{1}{N-1}\left\langle \overline{Df},\Delta^{(\epsilon)}\overline{x}\right\rangle\eqno{(144.1)}
$$
and the equation
\begin{equation}
\left\langle\overline{Df},\Delta^{(\epsilon)}\overline{x}\right\rangle=0\textrm{ is equivalent to  }\Delta^{(\epsilon)}f=0.
\end{equation}
In particular if $\overline{Df}=\overline{0}$, then $\Delta^{(\epsilon)}f=0$.\\
\\
\textbf{Theorem 27.1}\\
If $\lambda_{ij}$ is antisymmetric and 
\begin{equation}
\Delta^{(\lambda)}\overline{x}=\sum^{N}_{k=1}\sigma^{(\lambda)}_k\overline{e}_{k},
\end{equation}
then
\begin{equation}
\Delta^{(\lambda)}f=\sum^{N}_{k=1}\sigma^{(\lambda)}_{k}\nabla_{k}f
\end{equation}
\\
\textbf{Lemma 2.}\\
If $\lambda_{ij}$ is symmetric, then
\begin{equation}
\Delta^{(\lambda)}f=(N-1)\sum_{i\leq j}^{*}\lambda_{ij}\nabla_i\nabla_jf+\sum^{N}_{k=1}\left((N-1)A^{(\lambda)}_k-\frac{1}{N-1}\sum^{N}_{i=1}\lambda_{ik}R_i\right)(\nabla_kf),
\end{equation}
where the asterisk on the summation means that when $i<j$ the summands are multiplied with 2 and when $i=j$ with 1.\\
\\
\textbf{Proof.}\\
Assume that $\lambda_{ij}$ is symmetric, then we can write
$$
\Delta^{(\lambda)}f=\sum^{N}_{i,j=1}\lambda_{ij}D_i\nabla_jf=
$$
$$
=\sum^{N}_{k=1}\lambda_{kk}D_k\nabla_kf+\sum_{i<j}\lambda_{ij}\left(D_{i}\nabla_jf+D_j\nabla_if\right).
$$
But
$$
D_i\nabla_jf+D_j\nabla_if=(N-1)\nabla_i\nabla_jf-\frac{R_i\nabla_jf}{N-1}+(N-1)\nabla_j\nabla_if-\frac{R_j\nabla_if}{N-1}=
$$
$$
=(N-1)\left(\nabla_i\nabla_jf+\nabla_j\nabla_if\right)-\frac{1}{N-1}\left(R_i\nabla_jf+R_j\nabla_if\right)=
$$
$$
(N-1)\left(2\nabla_i\nabla_jf+\sum^{N}_{k=1}\nabla_kfQ_{ikj}\right)-\frac{1}{N-1}\left(R_i\nabla_jf+R_j\nabla_if\right).
$$
Hence
$$
\sum_{i<j}\lambda_{ij}\left(D_{i}\nabla_jf+D_j\nabla_if\right)=
$$
$$
=2(N-1)\sum_{i<j}\lambda_{ij}\nabla_i\nabla_jf+(N-1)\sum^{N}_{k=1}(\nabla_kf)\sum_{i<j}\lambda_{ij}Q_{ikj}-
$$
$$
-\frac{1}{N-1}\sum^{N}_{i,k=1}\lambda_{ik}R_i\nabla_kf+\frac{1}{N-1}\sum^{N}_{k=1}\lambda_{kk}R_k\nabla_kf=
$$
$$
=2(N-1)\sum_{i<j}\lambda_{ij}\nabla_i\nabla_jf+(N-1)\sum^{N}_{k=1}A^{(\lambda)}_{k}(\nabla_kf)-\frac{1}{N-1}\sum^{N}_{i,k=1}\lambda_{ik}R_i\nabla_kf+
$$
$$
+\frac{1}{N-1}\sum^{N}_{k=1}\lambda_{kk}R_k\nabla_kf.
$$
Also
$$
\sum^{N}_{k=1}\lambda_{kk}D_k\nabla_kf=(N-1)\sum^{N}_{k=1}\lambda_{kk}\nabla^2_kf-\frac{1}{N-1}\sum^{N}_{k=1}\lambda_{kk}R_k(\nabla_kf).
$$
Hence combining the above we get the first result.\\
\\
\textbf{Theorem 28.}\\
If $\lambda_{ij}$ is symmetric, then
\begin{equation}
\Delta^{(\lambda)}f=(N-1)\sum^{*}_{i\leq j}\lambda_{ij}\left(\nabla_jf\right)_{;i}+\left\langle \Delta^{(\lambda)}\overline{x},\overline{\textrm{grad}f}\right\rangle,
\end{equation}
where the asterisk in the sum means that if $i<j$, then the summands are multiplied with 2. In case $i=j$ with 1.\\ 
\\
\textbf{Proof.}\\
From Lemma we have that if $\lambda_{ij}$ is symmetric, then
$$
\Delta^{(\lambda)}f=(N-1)\sum_{i\leq j}^{*}\lambda_{ij}\nabla_i\nabla_jf+\sum^{N}_{k=1}\left((N-1)A^{(\lambda)}_k-\frac{1}{N-1}\sum^{N}_{i=1}\lambda_{ik}R_i\right)(\nabla_kf).
$$
But also from Theorem 22 we have
$$
\frac{1}{N-1}\sum^{N}_{i=1}\lambda_{ik}R_i=(N-1)\sum^{N}_{i,j=1}\lambda_{ij}q_{ikj}-\left\langle \Delta^{(\lambda)}\overline{x},\overline{e}_k\right\rangle
$$
Hence
$$
\sum^{N}_{k=1}\left((N-1)A^{(\lambda)}_k-\frac{1}{N-1}\sum^{N}_{i=1}\lambda_{ik}R_i\right)(\nabla_kf)=
$$
$$
=(N-1)\sum^{N}_{k=1}\sum_{i<j}\lambda_{ij}Q_{ikj}\nabla_kf-(N-1)\sum^{N}_{i,j,k=1}\lambda_{ij}q_{ikj}\nabla_kf+\sum^{N}_{k=1}\left\langle \Delta^{(\lambda)}\overline{x},\overline{e}_k\right\rangle\nabla_kf=
$$
$$
=(N-1)\sum^{N}_{k=1}\nabla_kf\left(\sum_{i<j}\lambda_{ij}q_{ikj}-\sum_{i<j}\lambda_{ij}q_{jki}-\sum_{i<j}\lambda_{ij}q_{ikj}-\sum_{i>j}\lambda_{ij}q_{ikj}\right)-
$$
$$
-(N-1)\sum^{N}_{k,i=1}\lambda_{ii}q_{iki}\nabla_kf+\left\langle \Delta^{(\lambda)}\overline{x},\overline{\textrm{grad}f}\right\rangle=
$$
$$
=-2(N-1)\sum^{N}_{k=1}\sum_{i<j}\lambda_{ij}q_{jki}\nabla_kf-(N-1)\sum^{N}_{k,i=1}\lambda_{ii}q_{iki}\nabla_kf+\left\langle \Delta^{(\lambda)}\overline{x},\overline{\textrm{grad}f}\right\rangle.
$$
Hence
$$
\Delta^{(\lambda)}f=2(N-1)\sum_{i< j}\lambda_{ij}\nabla_i\nabla_jf+(N-1)\sum^{N}_{i=1}\lambda_{ii}\nabla_i^2f
-2(N-1)\sum_{i<j}\sum^{N}_{k=1}\lambda_{ij}q_{jki}\nabla_kf-
$$
$$
-(N-1)\sum^{N}_{k,i=1}\lambda_{ii}q_{iki}\nabla_kf+\left\langle \Delta^{(\lambda)}\overline{x},\overline{\textrm{grad}f}\right\rangle=
$$
$$
=(N-1)\sum^{*}_{i\leq j}\lambda_{ij}\left(\nabla_jf\right)_{;i}+\left\langle \Delta^{(\lambda)}\overline{x},\overline{\textrm{grad}f}\right\rangle.
$$
\\
\textbf{Definition 11.}\\
Assume that $g_{ij}$ is the metric tensor of the surface $\textbf{S}$. Then  
\begin{equation}
\Delta^{(g)}\overline{x}=\overline{0}
\end{equation}
Iff
\begin{equation} (N-1)g_{ij}Q^{*}_{ikj}-\frac{1}{N-1}\sum^{N}_{l=1}\epsilon^{(g)}_{kl}Q_{ilj}=0,
\end{equation}
where
\begin{equation}
\epsilon^{(g)}_{kl}:=\sum^{N}_{i=1}g_{ik}\epsilon_{il}.
\end{equation}
We call a space $\textbf{S}$, $G-$space iff $\Delta^{(g)}(\overline{x})=\overline{0}$.\\
\\
\textbf{Theorem 29.}\\
If $\textbf{S}$ is a $G-$space, then for all functions $f$, we have
\begin{equation}
\Delta^{(g)}f=(N-1)\sum^{*}_{i\leq j}g_{ij}\left(\nabla_jf\right)_{;i}
\end{equation}
\\
\textbf{Theorem 30.}\\
If $\lambda_{ij}$ is any antisymmetric field, then
\begin{equation}
\sum^{N}_{i,j=1}\lambda_{ij}[D_i,\nabla_j]f+\frac{1}{f}\sum^{N}_{k=1}A^{(\lambda)}_kD_k\left(f^2\right)
=\frac{f}{2(N-1)}\sum^{N}_{i,j=1}\lambda_{ij}\Theta^{(2)}_{ij}(R)
.
\end{equation}
In case that $\textbf{S}$ is $\textbf{S}$-$R$, then
\begin{equation}
\sum^{N}_{i,j=1}\lambda_{ij}[D_i,\nabla_j]f=-\frac{1}{f}\sum^{N}_{k=1}A_k^{(\lambda)}D_k\left(f^2\right).
\end{equation}
\\
\textbf{Proof.}\\
We first evaluate the bracket.
$$
[D_i,\nabla_j]f=D_i\nabla_jf-\nabla_jD_if=(N-1)\nabla_i\nabla_jf-\frac{1}{N-1}R_i\nabla_jf-
$$  
$$
-\left((N-1)\nabla_j\nabla_if-\frac{1}{N-1}\nabla_j(R_if)\right)=
$$
$$
(N-1)\left(\nabla_i\nabla_jf-\nabla_j\nabla_if\right)-\frac{1}{N-1}R_i\nabla_jf+\frac{1}{N-1}(f\nabla_jR_i+R_i\nabla_jf)=
$$
$$
=-(N-1)\sum^{N}_{k=1}(\nabla_kf)Q_{ikj}+\frac{f}{N-1}\nabla_jR_i.
$$
Hence if we multiply with $\lambda_{ij}$ and sum with respect to $i,j$, we get
$$
P=\sum^{N}_{i,j=1}\lambda_{ij}[D_i,\nabla_j]f=-(N-1)\sum^{N}_{k,i,j=1}\lambda_{ij}\nabla_kfQ_{ikj}+\frac{f}{N-1}\sum^{N}_{i,j=1}\lambda_{ij}\nabla_jR_i.
$$
Also we have from Theorem 7 and the antisymmetric property of $\lambda_{ij}$:
$$
\sum^{N}_{i,j=1}\lambda_{ij}\nabla_jR_i=\frac{1}{2}\sum^{N}_{i,j=1}\lambda_{ij}(\nabla_jR_i-\nabla_iR_j)=-\frac{1}{2}\sum^{N}_{i,j,k=1}\lambda_{ij}R_kQ_{jki}+
$$
$$
+\frac{1}{2}\sum^{N}_{i,j=1}\lambda_{ij}\Theta^{(2)}_{ij}(R)
=\frac{1}{2}\sum^{N}_{i,j,k=1}\lambda_{ij}R_kQ_{ikj}+\frac{1}{2}\sum^{N}_{i,j=1}\lambda_{ij}\Theta^{(2)}_{ij}(R).
$$
Hence combining the above two results, we get
$$
P=-(N-1)\sum^{N}_{k,i,j=1}\lambda_{ij}\nabla_kfQ_{ikj}+\frac{f}{2(N-1)}\sum^{N}_{i,j,k=1}\lambda_{ij}R_kQ_{ikj}+
$$
$$
+\frac{f}{2(N-1)}\sum^{N}_{i,j=1}\lambda_{ij}\Theta^{(2)}_{ij}(R)\Rightarrow
$$
$$
fP=-\frac{1}{2}\sum^{N}_{k,i,j=1}\lambda_{ij}Q_{ikj}\left((N-1)\nabla_k\left(f^2\right)-\frac{R_kf^2}{N-1}\right)+
$$
$$
+\frac{f^2}{2(N-1)}\sum^{N}_{i,j=1}\lambda_{ij}\Theta_{ij}^{(2)}(R)
$$
and the result follows.\\
\\
\textbf{Corollary 11.}\\
If $\lambda_{ij}$ is antisymmetric, then
\begin{equation}
\sum_{i<j}\lambda_{ij}[D_i,\nabla_j]f=-(N-1)\sum^{N}_{k=1}A_k^{(\lambda)}\nabla_kf+\frac{f}{N-1}\sum_{i<j}\lambda_{ij}\nabla_jR_i.
\end{equation}
\\
\textbf{Theorem 31.}\\
If $\lambda_{ij}=\lambda^{+}_{ij}+\lambda^{-}_{ij}$, is any field, then
$$
\Pi^{(\lambda)} f:=\sum^{N}_{i,j=1}\lambda_{ij}\left[D_i,\nabla_j\right]f=\frac{1}{f(N-1)}\sum^{N}_{k=1}H^{(\lambda^{-})}_kD_k\left(f^2\right)+
$$
\begin{equation}
+\frac{f}{2(N-1)}\sum^{N}_{i,j=1}\lambda^{-}_{ij}\Theta^{(2)}_{ij}(R)
+\frac{f}{N-1}\sum^{N}_{i,j=1}\lambda^{+}_{ij}\nabla_jR_i.
\end{equation}
\\
\textbf{Corollary 12.}\\
If $\lambda_{ij}$ is any field and $\overline{D\left(f^2\right)}=\overline{0}$, then exists $\mu$ independed of $f$ such
\begin{equation}
\Pi^{(\lambda)}f=\mu f.
\end{equation}
In particular
\begin{equation}
\mu=\frac{1}{2(N-1)}\sum^{N}_{i,j=1}\lambda^{-}_{ij}\Theta^{(2)}_{ij}(R)
+\frac{1}{N-1}\sum^{N}_{i,j=1}\lambda^{+}_{ij}\nabla_jR_i.
\end{equation}
\\
\textbf{Theorem 32.}\\
If $\lambda_{ij}$ is symmetric, then $\Pi^{(\lambda)}$ is simplified considerably
\begin{equation}
\Pi^{(\lambda)}f=\frac{f}{N-1}\sum^{N}_{i,j=1}\lambda_{ij}\nabla_jR_i.
\end{equation}
\\
\textbf{Theorem 33.}\\
If $\lambda_{ij}$ is any symmetric field, then for every function $f$ we have
\begin{equation}
\Pi^{(\lambda)}f=0\textrm{, }\forall f
\end{equation}
iff
\begin{equation}
\sum^{N}_{i,j=1}\lambda_{ij}\nabla_jR_i=0.
\end{equation}
\\
\textbf{Remark 6.}\\
1) In any space $\textbf{S}$ we define the quantity   
\begin{equation}
\eta^{(\lambda)}:=\frac{1}{N-1}\sum^{N}_{i,j=1}\lambda_{ij}\nabla_jR_i.
\end{equation}
2) If $\lambda_{ij}$ is symmetric, then
\begin{equation}
\Pi^{(\lambda)}f=\eta^{(\lambda)}f.
\end{equation}
3) If $\lambda_{ij}=g_{ij}$ (the metric tensor), then $\eta^{(g)}:=\eta$, $\Pi^{(g)}:=\Pi$ and 
\begin{equation}
\Pi f:=\sum^{N}_{i,j=1}g_{ij}[D_i,\nabla_j]f=\eta f,
\end{equation}
where 
\begin{equation}
\eta=\frac{1}{N-1}\sum^{N}_{i,j=1}g_{ij}\nabla_jR_i.
\end{equation}
\\
\textbf{Theorem 34.}\\
If $\lambda_{ij}$ is any field and $\lambda^{(S)}_{ij}=\frac{1}{2}\left(\lambda_{ij}+\lambda_{ji}\right)$, then
\begin{equation}
\Delta^{(\lambda)}f=(N-1)\sum^{*}_{i\leq j}\lambda^{(S)}_{ij}\left(\nabla_jf\right)_{;i}+\left\langle\Delta^{(\lambda)}\overline{x},\overline{\textrm{grad}f}\right\rangle.
\end{equation}
\\
\textbf{Proof.}\\
Write $\lambda_{ij}=\lambda^{(S)}_{ij}+\lambda^{(A)}_{ij}$, where $\lambda^{(S)}_{ij}=\frac{1}{2}(\lambda_{ij}+\lambda_{ji})$ is the symmetric part of $\lambda_{ij}$ and $\lambda^{(A)}_{ij}=\frac{1}{2}(\lambda_{ij}-\lambda_{ji})$ is the antisymmetric part of $\lambda_{ij}$. Then we have
$$
\Delta^{(\lambda)}f=\sum^{N}_{i,j=1}\lambda^{(S)}_{ij}D_i\nabla_{j}f+\sum^{N}_{i,j=1}\lambda^{(A)}_{ij}D_i\nabla_{j}f=
$$ 
$$
=2(N-1)\sum_{i<j}\lambda^{(S)}_{ij}\nabla_i\nabla_jf+(N-1)\sum^{N}_{i=1}\lambda^{(S)}_{ii}\left(\nabla_if\right)_{;i}-2(N-1)\sum^{N}_{k=1}\sum_{i<j}\lambda^{(S)}_{ij}q_{jki}\nabla_kf+
$$
$$
+(N-1)\sum^{N}_{k,i,j=1}\lambda^{(S)}_{ij}q_{ikj}\nabla_kf-\frac{1}{N-1}\sum^{N}_{i,k=1}\lambda^{(S)}_{ik}R_i\nabla_kf-
$$
$$
-(N-1)\sum^{N}_{k=1}A^{(A)}_k\nabla_kf
-\frac{1}{N-1}\sum^{N}_{i,k=1}\lambda^{(A)}_{ik}R_i\nabla_kf=
$$
$$
=2(N-1)\sum_{i<j}\lambda^{(S)}_{ij}\nabla_i\nabla_jf+(N-1)\sum^{N}_{i=1}\lambda^{(S)}_{ii}\left(\nabla_if\right)_{;i}-
$$
$$
-(N-1)\sum^{N}_{k=1}\sum_{i>j}\lambda_{ij}\left(q_{ikj}-q_{jki}\right)\nabla_{k}f+(N-1)\sum^{N}_{k,i=1}\lambda_{ii}q_{iki}\nabla_kf-
$$
$$
-\frac{1}{N-1}\sum^{N}_{k,i=1}\lambda_{ik}R_i\nabla_kf=
$$
$$
=2(N-1)\sum_{i<j}\lambda^{(S)}_{ij}\nabla_i\nabla_jf+(N-1)\sum^{N}_{i=1}\lambda^{(S)}_{ii}\left(\nabla_if\right)_{;i}-
$$
$$
-(N-1)\sum^{N}_{k=1}\sum_{i>j}\lambda_{ij}\left(q_{ikj}-q_{jki}\right)\nabla_{k}f+(N-1)\sum^{N}_{k,i=1}\lambda_{ii}q_{iki}\nabla_kf+
$$
$$
+\left\langle \Delta^{(\lambda)}\overline{x},\overline{\textrm{grad}f}\right\rangle-(N-1)\sum^{N}_{k,i,j=1}\lambda_{ij}q_{jki}\nabla_kf=
$$
$$
=2(N-1)\sum_{i<j}\lambda^{(S)}_{ij}\nabla_i\nabla_jf+(N-1)\sum^{N}_{i=1}\lambda^{(S)}_{ii}\left(\nabla_if\right)_{;i}-
$$
$$
-(N-1)\sum^{N}_{k=1}\sum_{i>j}\lambda_{ij}\left(q_{ikj}-q_{jki}\right)\nabla_{k}f+\left\langle \Delta^{(\lambda)}\overline{x},\overline{\textrm{grad}f}\right\rangle-
$$
$$
-(N-1)\sum^{N}_{k=1}\sum_{i<j}\lambda_{ij}q_{jki}\nabla_kf-(N-1)\sum^{N}_{k=1}\sum_{i>j}\lambda_{ij}q_{jki}\nabla_kf=
$$
$$
=2(N-1)\sum_{i<j}\lambda^{(S)}_{ij}\nabla_i\nabla_jf+(N-1)\sum^{N}_{i=1}\lambda^{(S)}_{ii}\left(\nabla_if\right)_{;i}-
$$
$$
-(N-1)\sum^{N}_{k=1}\sum_{i<j}\lambda_{ji}q_{jki}\nabla_{k}f+\left\langle \Delta^{(\lambda)}\overline{x},\overline{\textrm{grad}f}\right\rangle-
$$
$$
-(N-1)\sum^{N}_{k=1}\sum_{i<j}\lambda_{ij}q_{jki}\nabla_kf=
$$
$$
=2(N-1)\sum_{i<j}\lambda^{(S)}_{ij}\nabla_i\nabla_jf+(N-1)\sum^{N}_{i=1}\lambda^{(S)}_{ii}\left(\nabla_if\right)_{;i}-
$$
$$
-2(N-1)\sum^{N}_{k=1}\sum_{i<j}\lambda^{(S)}_{ij}q_{jki}\nabla_{k}f+\left\langle \Delta^{(\lambda)}\overline{x},\overline{\textrm{grad}f}\right\rangle=
$$
$$
=(N-1)\sum^{*}_{i\leq j}\lambda^{(S)}_{ij}\left(\nabla_jf\right)_{;i}+\left\langle \Delta^{(\lambda)}\overline{x},\overline{\textrm{grad}f}\right\rangle.
$$
\\
\textbf{Definition 12.}\\
We call mean curvature vector of the surface $\textbf{S}$ the vector
\begin{equation}
\overline{H}=\sum^{N}_{k=1}H_k\overline{e}_k,
\end{equation}
where
\begin{equation}
H_k=(N-1)\sum^{N}_{i,j=1}g_{ij}q_{ikj}.
\end{equation}
\\
\textbf{Theorem 35.}\\
In case $\lambda_{ij}=g_{ij}$, then we have
\begin{equation}
\Delta^{(g)} f=(N-1)\sum^{N}_{i,j=1}g_{ij}\nabla_i\nabla_jf+\left\langle \Delta^{(g)}\overline{x}-\overline{H},\overline{\textrm{grad}f}\right\rangle.
\end{equation}
\\
\textbf{Proof.}\\
We know that $(\nabla_jf)_{;i}=(\nabla_if)_{;j}$ and $\lambda_{ij}=g_{ij}$ is symmetric. Hence from Theorem 32 we have
$$
\Delta^{(g)}f=2(N-1)\sum_{i<j}g_{ij}(\nabla_jf)_{;i}+(N-1)\sum^{N}_{i=1}g_{ii}(\nabla_if)_{;i}+\left\langle \Delta^{(g)}\overline{x},\overline{\textrm{grad}f}\right\rangle=
$$
$$
=(N-1)\sum^{N}_{i<j}g_{ij}(\nabla_jf)_{;i}+(N-1)\sum^{N}_{i>j}g_{ij}(\nabla_jf)_{;i}+(N-1)\sum^{N}_{i=1}g_{ii}(\nabla_if)_{;i}+
$$
$$
+\left\langle \Delta^{(g)}\overline{x},\overline{\textrm{grad}f}\right\rangle=
$$
$$
=(N-1)\sum^{N}_{i,j=1}g_{ij}(\nabla_jf)_{;i}+\left\langle \Delta^{(g)}\overline{x},\overline{\textrm{grad}f}\right\rangle=
$$
$$
=(N-1)\sum^{N}_{i,j=1}g_{ij}\nabla_i\nabla_jf-(N-1)\sum^{N}_{i,j=1}g_{ij}\sum^{N}_{k=1}q_{jki}\nabla_kf+\left\langle \Delta^{(g)}\overline{x},\overline{\textrm{grad}f}\right\rangle=
$$
$$
=(N-1)\sum^{N}_{i,j=1}g_{ij}\nabla_i\nabla_jf-\sum^{N}_{k=1}H_k\nabla_kf+\left\langle \Delta^{(g)}\overline{x},\overline{\textrm{grad}f}\right\rangle=
$$
$$
=(N-1)\sum^{N}_{i,j=1}g_{ij}\nabla_i\nabla_jf+\left\langle \Delta^{(g)}\overline{x}-\overline{H},\overline{\textrm{grad}f}\right\rangle
$$
and the theorem is proved.\\
\\
\textbf{Note 7.}\\
\textbf{i)} Actually the above theorem can be generalized for any $\lambda_{ij}$ as
\begin{equation}
\Delta^{(\lambda)}f=(N-1)\sum^{N}_{i,j=1}\lambda^{(S)}_{ij}\nabla_i\nabla_jf+\left\langle \Delta^{(\lambda)}\overline{x}-\overline{H}^{(S)},\overline{\textrm{grad}f}\right\rangle,
\end{equation}
where $\lambda^{(S)}_{ij}=\frac{\lambda_{ij}+\lambda_{ji}}{2}$, \begin{equation}
H^{(S)}_{k}=(N-1)\sum^{N}_{i,j=1}\lambda^{(S)}_{ij}q_{ikj}
\end{equation} 
and 
\begin{equation}
\overline{H}^{(S)}=\sum^{N}_{k=1}H^{(S)}_k\overline{e}_k.
\end{equation}
\textbf{ii)} If $\lambda_{ij}$ is symmetric, then
$$
\Delta^{(\lambda)}f=(N-1)\sum^{N}_{i,j=1}\lambda_{ij}\nabla_i\nabla_jf-\frac{1}{N-1}\sum^{N}_{k=1}R^{\{\lambda\}}_k\nabla_kf.\eqno{(173.1)}
$$
\textbf{iii)} If in a space $\textbf{S}$ we have for $\lambda_{ij}$ symmetric $\Delta^{(\lambda)}\left(\overline{x}\right)=\overline{H}^{(\lambda)}\Leftrightarrow\overline{R}^{\{\lambda\}}=\overline{0}$, then 
$$
\Delta^{(\lambda)}f=(N-1)\sum^{N}_{i,j=1}\lambda_{ij}\nabla_i\nabla_jf.
$$ 
\\
\textbf{Definition 13.}\\ We can also expand the definition of Beltrami operator acting to vectors. This can be done as follows:\\
If $\overline{A}=A_1\overline{\epsilon}_1+A_2\overline{\epsilon}_2+\ldots+A_{N}\overline{\epsilon}_{N}$, where $\{\overline{\epsilon}_i\}_{i=1,2,\ldots,N}$ is ''constant'' orthonormal base of $\textbf{E}=\textbf{R}^{N}$, then
\begin{equation}
\Delta_2(\overline{A})=\Delta_2(A_1)\overline{\epsilon}_1+\Delta_2(A_2)\overline{\epsilon}_2+\ldots+\Delta_2(A_{N})\overline{\epsilon}_{N}.
\end{equation}
\\
\textbf{Theorem 36.}\\
If $\overline{V}_1$ and $\overline{V}_2$ are vector fields, then 
\begin{equation}
\Delta_2\left\langle \overline{V}_1,\overline{V}_2\right\rangle=\left\langle \Delta_2\overline{V}_1,\overline{V}_2\right\rangle+\left\langle \overline{V}_1,\Delta_2\overline{V}_2\right\rangle+2(N-1)\sum^{N}_{k=1}\left\langle \nabla_k\overline{V}_1,\nabla_k\overline{V}_2\right\rangle.
\end{equation}
\\
\textbf{Proof.}\\
The proof follows by direct use of Theorem 7 and the identity 
\begin{equation}
\nabla_k\left\langle \overline{V}_1, \overline{V}_2\right\rangle=\left\langle\nabla_k\overline{V}_1,\overline{V}_2\right\rangle+
\left\langle \overline{V}_1,\nabla_k\overline{V}_2\right\rangle.
\end{equation}
\\
\textbf{Corollary 13.}
\begin{equation}
\left\langle \Delta_2\left(\overline{e}_k\right),\overline{e}_k \right\rangle=-(N-1)\sum^{N}_{j,l=1}q_{kjl}^2=\textrm{invariant.}
\end{equation}
\\
\textbf{Proof.}\\
From Theorem 36 and
$$
\left\langle \overline{e}_k,\overline{e}_k \right\rangle=1,
$$
we have
\begin{equation}
2 \left\langle \Delta_2\left(\overline{e}_k\right),\overline{e}_k \right\rangle+2(N-1)\sum^{N}_{l=1}\left|\nabla_l\left(\overline{e}_k\right)\right|^2=0.
\end{equation}
Hence we get the result.\\

From Theorem 15 and Proposition 1 one can easily see that
\begin{equation}
\Delta_2(\overline{x})=\sum^{N}_{k=1}\left((N-1)\sum^{N}_{l=1}q_{lkl}-\frac{1}{N-1}R_k\right)\overline{e}_k.
\end{equation}
Hence if we define as ''2-mean curvature vector'' the (not to confused with $h_k$ of ($eq1$),($eq2$) in pg.17):
\begin{equation}
\overline{h}=\sum^{N}_{k=1}h_k\overline{e}_k,
\end{equation}
where
\begin{equation}
h_k:=(N-1)\sum^{N}_{l=1}q_{lkl},
\end{equation}
then
\begin{equation}
h_k^{*}=\left\langle \Delta_2(\overline{x}),\overline{e}_k\right \rangle=h_k-\frac{1}{N-1}R_k
\end{equation}
and
\begin{equation}
\Delta_2\left(\overline{x}\right)=\sum^{N}_{k=1}h^{*}_k\overline{e}_k.
\end{equation}
Also we can write
$$
\Delta_2\overline{x}=\overline{h}-\frac{1}{N-1}\overline{R}
$$
and
$$
\left\langle\Delta_2\overline{x},\overline{\textrm{grad}(f)}\right\rangle=\left\langle\overline{H},\overline{\textrm{grad}(f)}\right\rangle+\Delta_2f-(N-1)\sum^{N}_{k=1}\nabla_k^2f.
$$
\\
\textbf{Proposition 6.}\\
The quantities $h_k=(N-1)\sum^{N}_{i=1}q_{iki}$ and $h_{ij}=\sum^{N}_{k=1}q_{ikj}^2$ are invariants.\\
\\
\textbf{Proof.}\\ It follows from invariant property of $\left\langle \nabla_l \overline{Y},\overline{e}_k\right \rangle$ for all $\overline{Y}$.\\
\\
\textbf{Proposition 7.}\\
For the mean curvature holds the following relation
\begin{equation}
h_k=\frac{1}{N-1}R_{k}+\left\langle\Delta_2(\overline{x}),\overline{e}_{k}\right\rangle=\textrm{invariant.}
\end{equation}
\\
\textbf{Definition 14.}\\
We call $\textbf{S}$ 2-minimal if $h_k=0$, $\forall k=1,2,\ldots,N$.\\ 
\\
\textbf{Theorem 37.}\\
The invariants $R^{\{M\}}_{j}=(N-1)\sum^{N}_{i=1}\left(\kappa^{\{M\}}_{ij}\right)^2$ have interesting properties. The sum  $R^{\{M\}}_{o}=\sum^{N}_{j=1}R^{\{M\}}_j$ is also invariant and 
\begin{equation}
\left|\nabla_k(\overline{e}_M)\right|^2=\frac{R^{\{M\}}_k}{N-1}
\end{equation}
and
\begin{equation}
\left \langle \Delta_2(\overline{e}_{M}),\overline{e}_{M}\right \rangle=-R^{\{M\}}_{o}.
\end{equation}
\\
\textbf{Proof.}\\
We have 
$$
\left\langle\nabla_l\overline{e}_k,\nabla_l\overline{e}_k\right\rangle=\sum^{N}_{s=1}q_{ksl}^2,
$$
which gives (185).\\ 
From relations (177),(178) and  (185) we get (186).\\ 
\\ 
\textbf{Definition 15.}\\
We call the $R^{\{M\}}_{o}$ as $R^{\{M\}}_{o}-$curvature.\\
\\

The extremely case $R^{\{M\}}_{o}=0$ for a certain $M$ happens if and only if $\kappa^{\{M\}}_{ij}=0$, for all $i,j=1,2,\ldots,N$. This leads from relation (10) to $\omega_{iM}=0$, for all $i=1,2,\ldots,N$ which means that $III_{M}=0$ and hence $\overline{e}_M$ is constant vector. We call such space flat in the $M$ direction.\\
\\
\textbf{Application 1.}\\
In the case $\overline{x}=\overline{e}_{N}$, then $\textbf{S}$ is a hypersphere and we have: 
\begin{equation}
\nabla_{k}\left(\overline{e}_N\right)=\overline{e}_k=\sum^{N}_{j=1}q_{Njk}\overline{e}_j.
\end{equation}
Hence
\begin{equation}
q_{mNk}=-\delta_{mk}.
\end{equation}
From Theorem 36 we have
$$
0=\Delta_2\left \langle \overline{e}_N,\overline{e}_N\right\rangle=2\left\langle \Delta_2\overline{e}_N,\overline{e}_N\right\rangle+2(N-1)N.
$$
Hence
\begin{equation}
R^{\{N\}}_{o}=-\left\langle \Delta_2\overline{e}_N,\overline{e}_N\right\rangle=N(N-1).
\end{equation}
Also
$$
h_{N}=-(N-1)\sum^{N}_{l=1}\delta_{ll}=-N(N-1)
$$
\begin{equation}
h^{*}_{N}=-N(N-1)=h_N-\frac{1}{N-1}R_{N}\Rightarrow R_{N}=0
\end{equation}
and
\begin{equation}
K^{\{N\}}=(-1)^N.
\end{equation}
From 
$$
R_{iNml}=-\sum^{N}_{s=1}rot_{lm}\left(q_{isl}q_{Nsm}\right)=\sum^{N}_{s=1}rot_{lm}\left(q_{isl}q_{sNm}\right)=
$$
$$
=-\sum^{N}_{s=1}rot_{lm}\left(q_{isl}\delta_{sm}-q_{ism}\delta_{sl}\right)=-\sum^{N}_{s=1}q_{isl}\delta_{sm}+\sum^{N}_{s=1}q_{ism}\delta_{sl}=
$$
\begin{equation}
=-q_{iml}+q_{ilm}=q_{mil}-q_{lim}=Q_{mil}=-Q_{lim}.
\end{equation}
Hence we get
$$
R_{NNml}=0.
$$
\\
\textbf{Theorem 38.}\\
In general
\begin{equation}
\Delta_2 f=(N-1)\sum^{N}_{k=1}\nabla_k^2f-\frac{1}{N-1}\left\langle \overline{R},\overline{\textrm{grad}f}\right\rangle.
\end{equation}
1) If $\textbf{S}$ is 2-minimal, then
\begin{equation}
\Delta_2\overline{x}=-\frac{1}{N-1}\overline{R}.
\end{equation}
2) If $t_i$ is any vector field, then easily in general
\begin{equation}
\sum^{N}_{i=1}(t_{i})_{;i}=\sum^{N}_{k=1}\nabla_kt_k-\sum^{N}_{k=1}t_kh_k
\end{equation}
and if $\textbf{S}$ is 2-minimal, then
$$
\sum^{N}_{k=1}\left(t_k\right)_{;k}=\sum^{N}_{k=1}\nabla_kt_k.
$$
\\
\textbf{Note 8.}\\
Assume that (with Einstein's notation)
$$
\Delta_2\Phi:=g^{lm}\left(\frac{\partial^2\Phi}{\partial x^l\partial x^m}-\Gamma^{a}_{lm}\frac{\partial \Phi}{\partial x^a}\right).\eqno{(a)}
$$
Then
$$
h_k=\left \langle \Delta_2(\overline{x}),\overline{e}_k\right\rangle=g^{lm}b_{kl,m}.\eqno{(b)}
$$
Hence $h_k$ is invariant. We call $h_k$ mean curvature tensor.\\
\\
\textbf{Proof.}\\
We know that
$$
b_{kl}=\left\langle \frac{\partial \overline{x}}{\partial x^l},\overline{e}_k\right\rangle\textrm{ and }
\partial_l\overline{e}_k=\Gamma^{a}_{lk}\overline{e}_{a}
$$

We differentiate the above first of two identities with respect to $x^m$ and we have
$$
\frac{\partial b_{kl}}{\partial x^m}=\left\langle \partial^2_{lm}\overline{x},\overline{e}_k\right\rangle+\left\langle\partial_{l}\overline{x},\partial_m\overline{e}_k\right\rangle=\left\langle \partial^2_{lm}\overline{x},\overline{e}_k\right\rangle+\Gamma^{a}_{mk}\left\langle\partial_l\overline{x},\overline{e}_a\right\rangle.
$$
Hence
$$
\left\langle \partial^2_{lm}\overline{x},\overline{e}_k\right\rangle=\partial_{m}b_{kl}-\Gamma^a_{mk}b_{al}
$$
But
$$
h_k=\left\langle \Delta_2(\overline{x}),\overline{e}_k\right\rangle=g^{lm}\left\langle\partial^2_{lm}\overline{x},\overline{e}_k\right\rangle-g^{lm}\Gamma^a_{lm}\left\langle\partial_a\overline{x},\overline{e}_k\right\rangle=
$$
$$
g^{lm}\left(\partial_mb_{kl}-\Gamma^a_{mk}b_{al}\right)-g^{lm}\Gamma^a_{lm}b_{ka}=
$$
$$
=g^{lm}\left(\partial_mb_{kl}-\Gamma^a_{km}b_{al}-\Gamma^a_{lm}b_{ka}\right)=g^{lm}b_{kl,m}.
$$

\section{The Spherical Forms}

Consider also the Pfaff derivatives of a function $f$ with respect to the form $\omega_{Mm}$. It holds
\begin{equation}
df=\sum^{N}_{k=1}(\widetilde{\nabla}_{k}f)\omega_{Mk}.
\end{equation}
Assume the connections $q^{(1)}_{mMj}$ such that 
\begin{equation}
\omega_{m}=\sum^{N}_{j=1}q^{(1)}_{mMj}\omega_{Mj}.
\end{equation}
Then from (10) we have
\begin{equation}
\omega_m=\sum^{N}_{s,j=1}q^{(1)}_{mMj}q_{Mjs}\omega_{s}.
\end{equation}
Hence
\begin{equation}
\sum^{N}_{j=1}q^{(1)}_{mMj}q_{Mjs}=\delta_{ms},
\end{equation}
where $\delta_{ij}$ is the usual Kronecker symbol. Using this (196) becomes
\begin{equation}
df=\sum^{N}_{k=1}\widetilde{\nabla}_kf\sum^{N}_{s=1}q_{Mks}\omega_s=\sum^{N}_{s=1}\left(\sum^{N}_{k=1}(\widetilde{\nabla}_kf)q_{Mks}\right)\omega_s.
\end{equation}
Hence
\begin{equation}
\nabla_sf=\sum^{N}_{k=1}(\widetilde{\nabla}_kf)q_{Mks}
\end{equation}
and using (199)
\begin{equation}
\widetilde{\nabla}_lf=\sum^{N}_{s=1}{\nabla}_sfq^{(1)}_{sMl}.
\end{equation}
We set $\widetilde{q}_{mjl}$ to be the connection
\begin{equation}
\widetilde{\nabla}_l(\overline{e}_m)=\sum^{N}_{j=1}\widetilde{q}_{mjl}\overline{e}_j.
\end{equation}
One can easily see that
\begin{equation}
\widetilde{q}_{mjl}=\sum^{N}_{s=1}q_{mjs}q^{(1)}_{sMl}.
\end{equation}
From $d\overline{e}_M=\sum^{N}_{k=1}\omega_{Mk} \overline{e}_k$ we get $\widetilde{\nabla}_k(\overline{e}_M)=\overline{e}_k$ and $\left\langle\widetilde{\nabla}^2_k(\overline{e}_M),\overline{e}_M\right\rangle=\widetilde{q}_{kMk}$\\

Also
if $w=\left\langle \overline{x},\overline{e}_M\right\rangle$, then
$$
dw=\left\langle d\overline{x},\overline{e}_M\right\rangle+\left\langle \overline{x},d\overline{e}_M\right\rangle=\omega_{M}+\left\langle \overline{x},\sum^{N}_{k=1}\widetilde{\nabla}_k(\overline{e}_M)\omega_{Mk}\right\rangle=
$$
$$
=\sum^{N}_{j=1}q^{(1)}_{MMj}\omega_{Mj}+\sum^{N}_{k=1}\left\langle \overline{x},\widetilde{\nabla}_k(\overline{e}_M)\right\rangle\omega_{Mk}=
$$
$$
=\sum^{N}_{j=1}q^{(1)}_{MMj}\omega_{Mj}+\sum^{N}_{k=1}\left\langle \overline{x},\overline{e}_k\right\rangle \omega_{Mk}.
$$
Hence
\begin{equation}
\widetilde{\nabla}_kw=q^{(1)}_{MMk}+\left\langle \overline{x},\overline{e}_k\right\rangle .
\end{equation}
Also 
$$
d\overline{x}=\sum^{N}_{k=1}(\nabla_k\overline{x})\omega_k=\sum^{N}_{k,j=1}\overline{e}_kq^{(1)}_{kMj}\omega_{Mj}=\sum^{N}_{j=1}\left(\sum^{N}_{k=1}q^{(1)}_{kMj}\overline{e}_k\right)\omega_{Mj}.
$$
From this we get
\begin{equation}
\widetilde{\nabla}_j\overline{x}=\sum^{N}_{k=1}q^{(1)}_{kMj}\overline{e}_k\textrm{, }q^{(1)}_{kMj}\textrm{ is invariant }
\end{equation}
and
$$
\widetilde{\nabla}_l\widetilde{\nabla}_l(w)=\widetilde{\nabla}_l\left(q^{(1)}_{MMl}\right)+\left\langle \widetilde{\nabla}_l\overline{x},\overline{e}_l\right\rangle+\left\langle \overline{x},\widetilde{\nabla}_l\overline{e}_l\right\rangle=
$$
$$
=\widetilde{\nabla}_l\left(q^{(1)}_{MMl}\right)+\left\langle \sum^{N}_{k=1}q^{(1)}_{kMl}\overline{e}_k,\overline{e}_l\right\rangle+\left\langle \overline{x},\sum^{N}_{j=1}\widetilde{q}_{ljl}\overline{e}_j\right\rangle=
$$
$$
=\widetilde{\nabla}_l\left(q^{(1)}_{MMl}\right)+q^{(1)}_{lMl}+\sum^{N}_{j=1}\left\langle\overline{x},\overline{e}_j\right\rangle\widetilde{q}_{ljl}
$$
But 
$$
\Delta_2^{III_{M}}w=(N-1)\sum^{N}_{j=1}\widetilde{\nabla}^2_jw-\frac{1}{N-1}\sum^{N}_{j=1}(\widetilde{\nabla}_jw)\widetilde{R}_j=
$$
\begin{equation}
=(N-1)\sum^{N}_{j=1}\widetilde{\nabla}_j\left(q^{(1)}_{MMj}\right)+(N-1)\sum^{N}_{l=1}q^{(1)}_{lMl}+(N-1)\sum^{N}_{l,k=1}\left\langle \overline{x},\overline{e}_k\right\rangle \widetilde{q}_{ljl}-
$$
$$
-\frac{1}{N-1}\sum^{N}_{j=1}\left\langle \overline{x},\overline{e}_j\right\rangle \widetilde{R}_j-\frac{1}{N-1}\sum^{N}_{j=1}q^{(1)}_{MMj}\widetilde{R}_j.
\end{equation}
Also
$$
\Delta_2^{III_{M}}\overline{e}_M=(N-1)\sum^{N}_{j=1}\widetilde{\nabla}^2_j\overline{e}_M-\frac{1}{N-1}\sum^{N}_{j=1}(\widetilde{\nabla}_j\overline{e}_M)\widetilde{R}_j=
$$
\begin{equation}
=(N-1)\sum^{N}_{l,j=1}\widetilde{q}_{ljl}\overline{e}_j-\frac{1}{N-1}\sum^{N}_{j=1}\widetilde{R}_j\overline{e}_j
\end{equation}
Hence
\begin{equation}
\left\langle\Delta_2^{III_{M}}\overline{e}_M,\overline{x}\right\rangle=(N-1)\sum^{N}_{l,j=1}\left\langle \overline{x},\overline{e}_j \right\rangle \widetilde{q}_{ljl}-\frac{1}{N-1}\sum^{N}_{j=1}\left\langle \overline{x},\overline{e}_j \right\rangle\widetilde{R}_{j}
\end{equation}
From (207) and (209) we get
$$
\Delta_2^{III_{M}}\left\langle \overline{x},\overline{e}_M\right\rangle-\left\langle \overline{x},\Delta_2^{III_{M}}\overline{e}_M\right\rangle=(N-1)\sum^{N}_{l=1}q^{(1)}_{lMl}+
$$
\begin{equation}
+(N-1)\sum^{N}_{j=1}\widetilde{\nabla}_j\left(q^{(1)}_{MMj}\right)-\frac{1}{N-1}\sum^{N}_{j=1}q^{(1)}_{MMj}\widetilde{R}_j.
\end{equation}
From Theorem 36 and formula (197) we get
$$
\Delta_2^{III_M}\left\langle \overline{x},\overline{e}_{M}\right\rangle=\left\langle \Delta_2^{III_M}\overline{x},\overline{e}_M\right\rangle+\left\langle \overline{x},\Delta^{III_M}_2\overline{e}_M\right\rangle+
$$
$$
+2(N-1)\sum^{N}_{k=1}\left\langle \widetilde{\nabla}_k\overline{x},\widetilde{\nabla}_k\overline{e}_M\right\rangle
$$
Or using (210)
$$
(N-1)\sum^{N}_{l=1}q^{(1)}_{lMl}+\sum^{N}_{j=1}\widetilde{D}_jq^{(1)}_{MMj}=\left\langle \Delta_2^{III_M}\overline{x},\overline{e}_M\right\rangle+
$$
$$
+2(N-1)\sum^{N}_{k=1}\left\langle \sum^{N}_{j=1}q^{(1)}_{jMk}\overline{e}_j,\overline{e}_k\right\rangle
$$
Or
$$
\sum^{N}_{j=1}\widetilde{D}_jq^{(1)}_{MMj}+(N-1)\sum^{N}_{l=1}q^{(1)}_{lMl}=\left\langle \Delta_2^{III_M}\overline{x},\overline{e}_M\right\rangle+
$$
$$
+2(N-1)\sum^{N}_{k=1} q^{(1)}_{kMk}.
$$
Hence we get the next\\
\\
\textbf{Theorem 39.}\\
In $\textbf{S}$ holds
\begin{equation}
\left\langle\Delta_2^{III_M}\overline{x},\overline{e}_M\right\rangle=-(N-1)\sum^{N}_{l=1}q^{(1)}_{lMl}+\sum^{N}_{l=1}\widetilde{D}_lq^{(1)}_{MMl}.
\end{equation}
In case $\omega_{M}=0$, then the above formula becomes
$$
\left\langle\Delta_2^{III_M}\overline{x},\overline{e}_M\right\rangle=-(N-1)\sum^{N}_{l=1}q^{(1)}_{lMl}.\eqno{(211.1)}
$$

\section{General Forms and Invariants}

Let $C$ be a curve (one dimensional object) of  $\textbf{S}$. Let also $s$ be the physical parameter of $C$. Then 
\begin{equation}
\overline{t}=\left(\frac{d\overline{x}}{ds}\right)_C,
\end{equation}
is the tangent vector of $C$ in $P\in \textbf{S}$. If we assume that $C$ lays in a hypersurface $S_{M-1}$ and we choose  $\overline{n}=\overline{e}_{M}$ to be the normal vector of the tangent space of $S_{M-1}$, then  
\begin{equation}
\left\langle \overline{t},\overline{n}\right\rangle=0.
\end{equation}
\\
\textbf{Definition 16.}\\
We call vertical curvature of $C\in S_{M-1}$ the quantity 
\begin{equation}
\left(\frac{1}{\rho^{*}}\right)_C=\left\langle \frac{d\overline{t}}{ds},\overline{n}\right\rangle.
\end{equation}
\\
\textbf{Theorem 40.}\\    
The vertical curvature $\left(\frac{1}{\rho^{*}}\right)_C$ is an invariant (by the word invariant we mean that remains unchanged in every acceptable transformation of the coordinates $u_i$).\\ 
\\
\textbf{Proof.}\\
This can be shown as follows. Derivate (213) to get
$$ 
\left\langle\frac{d\overline{t}}{ds},\overline{n}
\right\rangle+\left\langle \overline{t},\frac{d\overline{n}}{ds}\right\rangle=0.
$$
Hence
$$
\left(\frac{1}{\rho^{*}}\right)_{C}=-\left\langle \frac{d\overline{x}}{ds},\frac{d\overline{n}}{ds}\right\rangle=-\frac{II_{M}}{I}=-\frac{\sum^{N}_{k=1}\omega_k\omega_{Mk}}{\sum^{N}_{k=1}\omega_k^2}=\textrm{invariant}.
$$
\\

Using relations (8),(10), we get 
$$
\left(\frac{1}{\rho^{*}}\right)_{C_i}=-\left\langle \left(\frac{d\overline{x}}{ds}\right)_{C_i},\left(\frac{d\overline{e}_{M}}{ds}\right)_{C_i}\right\rangle=
$$
$$
=-\left\langle \sum^{N}_{k=1}\left(\frac{\omega_k}{ds}\right)_{C_i}\overline{e}_k,\sum^{N}_{l=1}\left(\frac{\omega_{Ml}}{ds}\right)_{C_i}\overline{e}_l\right\rangle
=-\sum^{N}_{k,l=1}\left(\frac{\omega_k}{ds}\right)_{C_i}\left(\frac{\omega_{Ml}}{ds}\right)_{C_l}\delta_{kl}=
$$
$$
=-\sum^{N}_{k=1}\left(\frac{\omega_k}{ds}\right)_{C_i}\left(\frac{\omega_{Mk}}{ds}\right)_{C_i}
=-\sum^{N}_{k,m=1}\left(\frac{\omega_k}{ds}\right)_{C_i}q_{Mkm}\left(\frac{\omega_m}{ds}\right)_{C_i}.
$$
Hence
\begin{equation}
\left(\frac{1}{\rho^{*}}\right)_{C_i}=\sum^{N}_{k,m=1}\kappa^{\{M\}}_{km}\left(\frac{\omega_k}{ds}\right)_{C_i}\left(\frac{\omega_m}{ds}\right)_{C_i}.
\end{equation}
Assume now, in general that exist $N-1$ curves $C_i$, $i=1,2\ldots,N-1$ passing through every point $P$ of $S_{M-1}$ and are vertical to each other. Then in all $P\in S_{M-1}$ we have
\begin{equation}
\left\langle \left(\frac{d\overline{x}}{ds}\right)_{C_i},\left(\frac{d\overline{x}}{ds}\right)_{C_j}\right\rangle=\delta_{ij}\textrm{, }i,j=1,2,\ldots,N-1
\end{equation} 
and
$$
\left\langle \left(\frac{d\overline{x}}{ds}\right)_{C_i},\overline{n} \right\rangle=0.
$$
Also for these curves hold  
\begin{equation}
\left(\frac{d\overline{x}}{ds}\right)_{C_i}=\sum^{N}_{k=1}\left(\frac{\omega_k}{ds}\right)_{C_i}\overline{e}_k=\sum^{N}_{k=1}\lambda_{ik}\overline{e}_{k},
\end{equation}
where we have set
\begin{equation}
\left(\frac{\omega_k}{ds}\right)_{C_i}=\lambda_{ik}.
\end{equation}
Clearly from the orthogonality of $C_i$ we have that (where we have assumed with no loss of generality that $\omega_{M}=0$ hence $\lambda_{iM}=0$):
$$
\sum_{1\leq k\leq N}^{*}\lambda_{ik}\lambda_{jk}=\delta_{ij}\textrm{, }\forall i,j\in\{1,2,\ldots,N-1\}.
$$  
Where the asterisc in the sumation means that the value $k=M$ is omited.\\
Assuming these facts in (215), we get (using the fact that any real unitary matrix is symmetric): 
$$
\sum_{1\leq i\leq N}^{*}\left(\frac{1}{\rho^{*}}\right)_{C_i}=\sum^{*}_{1\leq i\leq N}\sum^{N}_{k,m=1}\kappa^{\{M\}}_{km}\lambda_{ik}\lambda_{im}=\sum^{N}_{k,m=1}\kappa^{\{M\}}_{km}\delta_{km}=
$$
$$
=\sum^{N}_{k=1}\kappa^{\{M\}}_{kk}=\textrm{invariant}.
$$
From the above we get the next\\
\\
\textbf{Theorem 41.}\\
From any point $P$ in a hypersurface $S_{M-1}$ of the space $\textbf{S}$, there pass $N-1$ vertical curves $C_i,i=1,2,\ldots,N-1$ and their vertical curvatures satisfy 
\begin{equation}
\sum_{1\leq i\leq N}^{*}\left(\frac{1}{\rho^{*}}\right)_{C_i}=\frac{h_M}{N-1}.
\end{equation}     
\textbf{Remark.} We mention here that with the word invariant we mean any quantity that remains unchanged in every acceptable choose of parameters $\{u_i\}_{i=1,2,\ldots,N}$ as also change of position and rotation of $\textbf{S}$.\\
\\

Set 
\begin{equation}
T=\sum^{N}_{i=1}A_i\omega_{i}^2+\sum^{N}_{i,j=1}B_{ij}\omega_{ij}^2+\sum^{N}_{i,j,k=1}C_{ijk}\omega_{i}\omega_{jk}+\sum^{N}_{i,j,f,m=1}E_{ijfm}\omega_{ij}\omega_{fm}
\end{equation}
and assume the $N-1$ orthogonal curves $C_i$. If the direction of $C_i$ is
\begin{equation}
d_i=\left(\frac{T}{(ds)^2}\right)_{C_i}\textrm{, }i=1,2,\ldots,N-1\textrm{ and }d_{N}=0,
\end{equation} 
then summing in all $N$ directions we get, after using the orthogonality and making simplifications
$$
\sum^{N}_{i=1}d_i=\sum^{N}_{i=1}A_i+\sum^{N}_{i,j,m=1}B_{ij}q_{ijm}^2+\sum^{N}_{i,j,k=1}C_{ijk}q_{jki}+
$$
\begin{equation}
+\sum^{N}_{i,j,f,m,n=1}E_{ijfmn}q_{ijn}q_{fmn}
\end{equation}
As application we get that
\begin{equation}
T_0=AI+BII+CIII=A\sum^{N}_{i=1}\omega_{i}^2+B\sum^{N}_{i=1}\omega_{Ni}^2+C\sum^{N}_{i=1}\omega_{i}\omega_{Ni}=\textrm{invariant}.
\end{equation}
Then summing the vertical directions we get 
\begin{equation}
\sum^{N-1}_{i=1}d_i=A+B\sum^{N}_{i,j=1}q_{iNj}^2+C\sum^{N}_{i=1}q_{iNi}=\textrm{invariant}
\end{equation}
For every choice of costants $A,B,C$.\\
\\
\textbf{Definition 17.}\\
Let now $C$ be a space curve and $\overline{t}_i$, $i=1,2,\ldots,N$ is a family of orthonormal vectors of $C$. We define 
\begin{equation}
\left(\frac{1}{\rho_{ik}}\right)_{C}:=\left\langle\left( \frac{d\overline{t}_i}{ds}\right)_C,\overline{t}_k\right\rangle\textrm{, }\forall i,k\in\{1,2,\ldots,N\}
\end{equation}
and call
$\left(\frac{1}{\rho_{ik}}\right)_C$ as ''$ik-$curvature'' of $C$ and 
\begin{equation}
\frac{d\overline{t}_i}{ds}=\sum^{N}_{k=1}\left(\frac{1}{\rho_{ik}}\right)_C\overline{t}_k\textrm{, }\forall i=1,2,\ldots,N.
\end{equation} 
\\
\textbf{Theorem 42.}\\
The $\left
(\frac{1}{\rho_{ij}}\right)_C-$curvatures are 
semi-invariants of the space. Moreover        
\begin{equation}
\left(\frac{1}{\rho_{ij}}\right)_C=-\left(\frac{1}{\rho_{ji}}\right)_C
\end{equation}
and if
$$
\overline{t}_i=\sum^{N}_{k=1}A_{ki}\overline{e}_k,
$$
then
\begin{equation}
\left(\frac{1}{\rho_{ij}}\right)_C=\sum^{N}_{l=1}\sum^{N}_{m=1}A_{\{m\}i;l}A_{mj}\left(\frac{\omega_l}{ds}\right)_{C}, 
\end{equation}
where $A_{\{m\}i;l}=\nabla_lA_{mi}-\sum^{N}_{k=1}q_{mkl}A_{ki}$.\\
\\
\textbf{Proof.}\\
We express $\overline{t}_i$ in the base $\overline{e}_k$ and differentiate with respect to the canonical parameter $s$ of $C$, hence
$$
\overline{t}_i=\sum^{N}_{k=1}A_{ki}\overline{e}_k,
$$
where $A_{ki}$ is unitary matrix i.e.
$$
\sum^{N}_{i=1}A_{ki}A_{li}=\delta_{kl}.
$$
Hence
\begin{equation}
\overline{e}_j=\sum^{N}_{k=1}A^{T}_{kj}\overline{t}_k,
\end{equation}
where $A^{T}$ is the symmetric of $A$.\\ 
Differentiating with respect to $s$ of $\overline{t}_i$, we get
$$
\frac{d\overline{t}_i}{ds}=\sum^{N}_{k=1}\frac{dA_{ki}}{ds}\overline{e}_k+\sum^{N}_{k=1}A_{ki}\frac{d\overline{e}_k}{ds}=
$$
$$
=\sum^{N}_{k=1}\frac{dA_{ki}}{ds}\overline{e}_k+\sum^{N}_{k=1}A_{ki}\sum^{N}_{l=1}\nabla_l\overline{e}_k\frac{\omega_l}{ds}=
$$
$$
=\sum^{N}_{k=1}\frac{dA_{ki}}{ds}\overline{e}_k+\sum^{N}_{k,l,m=1}A_{ki}q_{kml}\frac{\omega_l}{ds}\overline{e}_m
=
$$
$$
=\sum^{N}_{k,l=1}\frac{(\nabla_lA_{ki})\omega_l}{ds}\overline{e}_k+\sum^{N}_{k,l,m=1}A_{ki}q_{kml}\frac{\omega_l}{ds}\overline{e}_m=
$$
$$
=\sum^{N}_{m,l=1}\frac{(\nabla_lA_{mi})\omega_l}{ds}\overline{e}_m+\sum^{N}_{k,l,m=1}A_{ki}q_{kml}\frac{\omega_l}{ds}\overline{e}_m.
$$ 
Hence we can write
\begin{equation}
c_{im}:=\left\langle\frac{d\overline{t}_i}{ds},\overline{e}_m\right\rangle=\sum^{N}_{l=1}\left(\nabla_lA_{mi}-\sum^{N}_{k=1}q_{mkl}A_{ki}\right)\frac{\omega_l}{ds}=\sum^{N}_{l=1}A_{\{m\}i;l}\frac{\omega_l}{ds}.
\end{equation}
Hence
\begin{equation}
c_{im}=\sum^{N}_{l=1}\nabla_lA_{mi}\left(\frac{\omega_l}{ds}\right)_C-\sum^{N}_{l,k=1}A_{ki}q_{lkm}\left(\frac{\omega_l}{ds}\right)_C.
\end{equation}
Hence
\begin{equation}
\frac{dA_{mi}}{ds}=\left\langle\frac{d\overline{t}_i}{ds},\overline{e}_m\right\rangle+\sum^{N}_{l,k=1}A_{ki}q_{lkm}\left(\frac{\omega_l}{ds}\right)_C
\end{equation}
Also
$$
\left\langle \overline{t}_i,\overline{t}_j\right\rangle=\delta_{ij}\Rightarrow \left\langle \frac{d\overline{t}_i}{ds},\overline{t}_j\right\rangle+\left\langle \overline{t}_i,\frac{d\overline{t}_j}{ds}\right\rangle=0.
$$
Hence
\begin{equation}
\left(\frac{1}{\rho_{ij}}\right)_{C}=-\left(\frac{1}{\rho_{ji}}\right)_{C}.
\end{equation}
From the orthonormality of $A_{ij}$ and $A^{T}_{ij}=A_{ji}$ we have
\begin{equation}
\sum^{N}_{i=1}A_{mi}A_{ji}=\delta_{mj}\Leftrightarrow \sum^{N}_{m=1}A_{im}A_{jm}=\delta_{ij}.
\end{equation}
Differentiating we get
\begin{equation}
\sum^{N}_{i=1}A_{\{m\}i;l}A_{ji}+\sum^{N}_{i=1}A_{mi}A_{\{j\}i;l}=0
\end{equation}
Moreover 
\begin{equation}
\frac{d\overline{t}_i}{ds}=\sum^{N}_{m=1}c_{im}\overline{e}_m=\sum^{N}_{j=1}\sum^{N}_{l,m=1}A_{\{m\}i;l}A_{mj}\left(\frac{\omega_l}{ds}\right)_C\overline{t}_j
\end{equation}
and the curvatures will be
\begin{equation}
\left(\frac{1}{\rho_{ij}}\right)_C=\sum^{N}_{l=1}\sum^{N}_{m=1}A_{\{m\}i;l}A_{mj}\left(\frac{\omega_{l}}{ds}\right)_{C}.
\end{equation}
Hence
$$
\sum^{N}_{j=1}A_{pj}\left(\frac{1}{\rho_{ij}}\right)_C=\sum^{N}_{l=1}\sum^{N}_{m,j=1}A_{\{m\}i;l}A_{mj}A_{pj}\left(\frac{\omega_{l}}{ds}\right)_{C}=
$$
$$
=\sum^{N}_{l=1}\sum^{N}_{m=1}A_{\{m\}i;l}\delta_{mp}\left(\frac{\omega_{l}}{ds}\right)_{C}=\sum^{N}_{l=1}A_{\{p\}i;l}\left(\frac{\omega_{l}}{ds}\right)_{C}=c_{ip}.
$$
This lead us to write
$$
\sum^{N}_{i,j=1}A_{ni}A_{pj}\left(\frac{1}{\rho_{ij}}\right)_C=\sum^{N}_{i=1}A_{ni}c_{ip}=\sum^{N}_{l,i=1}A_{\{p\}i;l}A_{ni}\left(\frac{\omega_l}{ds}\right)_{C}=
$$
$$
=\sum^{N}_{l,i=1}A^{T}_{\{i\}p;l}A^{T}_{in}\left(\frac{\omega_l}{ds}\right)_C=\sum^{N}_{l,m=1}A^{T}_{\{m\}p;l}A^{T}_{mn}\left(\frac{\omega_l}{ds}\right)_C=
$$
$$
=\left(\sum^{N}_{l,m=1}A_{\{m\}n;l}A_{mp}\left(\frac{\omega_l}{ds}\right)_C\right)^T=\left(\frac{1}{\rho_{np}}\right)^{T}_C=
$$
$$
=\left(\sum^{N}_{l,m=1}\left(\nabla_lA_{\{m\}n}A_{mp}-\sum^{N}_{k=1}q_{mkl}A_{kn}A_{mp}\right)\right)^T\left(\frac{\omega_l}{ds}\right)_C=
$$
$$
=\sum^{N}_{l,m=1}\left(\nabla_lA_{\{m\}p}A_{mn}-\sum^{N}_{k=1}q_{mkl}A_{kp}A_{mn}\right)\left(\frac{\omega_l}{ds}\right)_C=
$$
$$
=\sum^{N}_{l,m=1}A_{\{m\}p;l}A_{mn}\left(\frac{\omega_l}{ds}\right)_C=
$$
$$
=\left(\frac{1}{\rho_{pn}}\right)_C=-\left(\frac{1}{\rho_{np}}\right)_C.
$$
Since we have
$$
\left(A^{T}_{\{i\}j}\right)_{;l}=\nabla_l A_{ji}-\sum^{N}_{k=1}q_{jkl}A_{ki}
$$
and
$$
\left(A_{\{i\}j;l}\right)^{T}=\left(\nabla_lA_{ij}-\sum^{N}_{k=1}q_{ikl}A_{kj}\right)^{T}=\nabla_lA_{ji}-\sum^{N}_{k=1}q_{jkl}A_{ki}
$$
and
\begin{equation}
\left(A^{T}_{\{i\}j}\right)_{;l}=\left(A_{\{i\}j;l}\right)^{T}.
\end{equation}

\[
\]

\centerline{\bf References}

\[
\]

[1]: Nirmala Prakash. ''Differential Geometry An Integrated Approach''. Tata McGraw-Hill Publishing Company Limited. New Delhi. 1981.\\

[2]: Bo-Yu Hou, Bo-Yuan Hou. ''Differential Geometry for Physicists''. World Scientific. Singapore, New Jersey, London, Hong Kong. 1997.\\

[3]: N.K. Stephanidis. ''Differential Geometry''. Vol. I. Zitis Pub. Thessaloniki, Greece. 1995.\\

[4]: N.K. Stephanidis. ''Differential Geometry''. Vol. II. Zitis Pub. Thessaloniki, Greece. 1987.\\

[5]: D. Dimitropoulou-Psomopoulou. ''Calculus of Differential Forms''. 2nd edition. Zitis Pub. Thessaloniki, Greece. 1993.\\ 

[6]: E.A. Iliopoulou, F. Gouli-Andreou. ''Introduction to Riemann Geometry''. Zitis Pub. Thessaloniki, Greece. 1985.\\

[7]: N.K. Spyrou. ''Introduction to the General Theory of Ralativity''. Gartaganis Pub. Thessaloniki, Greece. 1989.\\

\end{document}